\documentclass[]{amsart}
\usepackage{graphicx}
\usepackage[T1]{fontenc}
\usepackage[latin1]{inputenc}
\usepackage{enumerate}
\usepackage{pstricks}
\usepackage{latexsym,amssymb,amsthm,amsxtra}
\newcommand\car{{\mathbf 1}}
\newcommand\N{{\mathbb N}}
\renewcommand\S{{\mathcal S}}
\newcommand\M{{\mathcal M}}
\newcommand\C{{\mathcal C}}
\newcommand\D{{\mathcal D}}
\renewcommand\P{{\mathbf P}} 
\newcommand{\B}{{\mathfrak B}}
\newcommand{\A}{{\mathcal  A}}

\newcommand{\R}{{\mathbb R}}

\newcommand\esp[1]{{\mathbf E}\left[#1\right]}
\newcommand\pr[1]{{\mathbf P}\left[#1\right]}
\newcommand{\thmref}[1]{Theorem~\ref{#1}}
\newcommand{\lemref}[1]{Lemma~\ref{#1}}

\newcommand{\req}[1]{(\ref{#1})}
\newcommand{\pae}[1]{\mbox{$\lfloor \kern-1pt #1 \kern-1pt \rfloor$}}
\newcommand{\paep}[1]{\mbox{$\lceil \kern-1pt #1 \kern-1pt \rceil$}}

\newcommand\F{{\mathcal F}}

\newcommand\ovnun{{\bar \nu^n}}
\newcommand\ovrt{{\bar r_t}} 

\newcommand\ovnunPD{{\bar \nu^{n,\text{\tiny{PD}}}}}

\newcommand\ovnuPD{{\bar \nu^{*,\text{\tiny{PD}}}}}

\newcommand\mFnt{{\mathcal F_t^n}}
\newcommand\mGnt{{\mathcal G_t^n}}

\newcommand\lan{{\lambda^n}}
\newcommand\mun{{\mu^n}}
\newcommand\ton{{\tau_0^n}}
\newcommand\ovton{{\bar \tau_0^n}}

\newcommand\omn{{\omega_0^{n}}}
\newcommand\ovomn{{\bar \omega_0^{n}}}

\newcommand\Rp{{\R^*_+}}

\newcommand\n{{^n}}
\newcommand\tend{{\underset{n \rightarrow \infty}{\longrightarrow}}}

\newcommand\cro[1]{\langle #1 \rangle}
\newcommand\proc[1]{\left(#1_t\right)_{t \ge 0}}
\newcommand\proce[1]{\left(#1\right)_{t \ge 0}}
\newcommand\procp[1]{\left(#1(t),t \ge 0\right)}
\newcommand\suite[1]{\left\{#1\right\}_{n \in \N^*}}
\newcommand\suitei[1]{\left\{#1\right\}_{i \in \N^*}}

\newcommand\barn[1]{\bar#1^{(n)}}
\newcommand\barnp[1]{\bar#1^{n}_{\phi}}
\renewcommand\A{{\mathcal A}}
\newtheorem{hypo}{Hypothesis}{\bf}{\it}
\newtheorem{theorem}{Theorem}
\newtheorem{definition}{Definition}
\newtheorem{lemma}{Lemma}
\newtheorem{proposition}{Proposition}
\newtheorem{corollary}{Corollary}
\begin{document}

\title[Fluid limit of a heavily loaded EDF queue]{Fluid limit of a heavily loaded EDF queue with impatient customers}

\author{L. Decreusefond}
\address{GET/ENST - CNRS UMR 5141 \\
46, rue Barrault\\
Paris, 75634, FRANCE}
\email{Laurent.Decreusefond@enst.fr}

\author{P. Moyal}
\address{GET/ENST - CNRS UMR 5141 \\
46, rue Barrault\\
Paris, 75634, FRANCE}
\email{Pascal.Moyal@enst.fr}


\begin{abstract}
  In this paper, we present the fluid limit of an heavily loaded
  {Earliest Deadline First} queue with impatient customers,
  represented by a measure-valued process keeping track of residual
  time-credits of lost and waiting customers. This fluid limit is the
  solution of an integrated transport equation.  We then use this fluid
  limit to derive fluid approximations of the processes counting
  the number of waiting and already lost customers.
\end{abstract}
\maketitle
\section{Introduction}
Queueing theory is a keystone of the development of current
telecommunications systems. Engineers now aim to guarantee the grade
of service customers are entitled to receive according to their
contracts with the carrier. One way to meet this objective is to
schedule requests according to their ``importance''. Of the utmost
interest, are the audio and video traffic flows, which are subject to
severe transmission delay constraints. In particular, some requests
can be thought as ``impatient'' since it may be better to discard some
packets which would eventually arrive too late in order to favor some
other packets which still can meet their delay requirements. 

 Another very active
branch of queueing theory is, nowadays, devoted to the analysis of call-centers
where customers are \textsl{impatient}: they tolerate to wait up to a
certain limit upon which they depart from the queueing line hence are
considered as lost both for the queueing system and for the
called-service provider. It is then crucial to develop service
disciplines which ensure a maximal number of served customers by
controlling waiting times keeping them within impatience bounds. These
disciplines are commonly referred to as real-time service disciplines.

In real-time queuing theory, each customer is not only
identified by his arrival time  and service duration but also by a
\textsl{deadline}. This means that a customer has a given period of
time (his {time credit}, {i.e.}, the remaining time before his 
deadline) during 
which he should enter the service booth. This time credit decreases at unit 
rate as time goes on. If it 
expires before the customer enters service, the customer is either lost
and the discipline is said to be hard, or he is kept in the
waiting line and the discipline is said to be soft. The discipline we address here is the so-called Earliest-Deadline-First discipline in its ``hard'' version:  the customer having the smallest
time credit is served first and  whenever the credit-time of a customer expires before it is served, this customer is lost.

To find which service discipline is best, one can compare them within
a static scenario, i.e., customers to be served are all present at
initial time and no new customer enter the system, service duration
and impatience of each customer are all known at the beginning; or in
dynamic environments, i.e., customers arrive randomly, their service
duration and impatience are only known stochastically. In both
settings, it appears that the so-called Earliest-Deadline-First (EDF
for short) discipline is optimal. It is known for a while~\cite{Dertou74}, that EDF discipline is
optimal for the static approach: if any (real-time) service discipline can serve the
customers of a given scenario without loss then EDF also does. Within
random environments, it has been proven in \cite{PanTow88} and
generalized in \cite{moyal05} that EDF discipline ensures the
least possible failure probability, i.e., the least number customers lost by missing their deadline.

Yet, apart from the notable exception of deterministic deadlines for
which EDF discipline reduces to the FIFO service policy with
impatience, no closed form of the loss probability is known.
The only satisfying quantitative approach so far consists
in numerically assessing the loss probability for an EDF system with Markov-chain approximations
\cite{Hon89,Pin91,Nain92,PanTow88}.

When no simple tractable object can describe a queueing system, one
wants to identify its ``mean behavior''. One hopes that a Markovian
process characterizing the system, when suitably normalized, can be
approximated by a {fluid limit} that is, a deterministic
continuous function of the time. Then the fluid limit describes the
general behavior of the considered process.  Numerous queueing systems
have already been investigated this way (see, for instance,
\cite{Rob,Bor67} for a pure delay system, \cite{Doytchi} for a
soft EDF queue, \cite{Gromoll} for a queue run under a processor
sharing service discipline).  For instance, $L_t$ being the amount of
customers at time $t$ in an $M/M/1$ queue with parameters $\lambda$
and $\mu$, one proves that the sequence of processes
$\suite{\bar{L}^n}$ defined by $\barn{L}_0=1$ and for all $t>0,$
$\barn{L}_t:=n^{-1}L_{nt}$ tends in distribution to
$\bar{L}:=\left((1+(\lambda -\mu)t)^+,\, t \in \R^+\right)$ and that
$\suite{\sqrt{n}\left(\barn{L}-\bar L\right)}$ converges in
distribution to a diffusion process. The fluid approximation of the
system presents the same first order characteristics as the ``real''
system: it fills in at velocity $\lambda$ and empties at velocity
$\mu$, the congestion reaches $0$ to the condition $\lambda<\mu$ (this
is Loynes's stability condition) after a time $\lambda - \mu$ (mean
duration of a busy period).

We want to obtain the same type of information for an $M/M/1$ queueing
system with impatient customers. In this case, it is easily seen 
that the process $\proc{X}$ which counts the number of customers in
the system is no longer Markovian. Indeed, the value of $X_{t+h}$ not
only depends on $X_t$, but also on all the time credits
 of all the $X_t$ customers
present in the queue at time $t$. Therefore we describe the system by
the point measure-valued process $\proc{\nu}$ whose unit of mass are
the time credits of all the customers waiting in the queue, or already
discarded.  

Formally, it is rather straightforward in our case to write down
the infinitesimal generator of the Markov process $\proc{\nu}$, see
\thmref{pro:Feller} below.  It is made of four terms, all but
one are standard and represent the evolution of the process when an
arrival, a departure or nothing occurs during an infinitesimal time
period. The natural but unusual term is the term due to the continuous
decreasing of the residual deadlines at unit rate as time goes
on. This term involves a ``spatial derivative'' of the measure $\nu$,
a notion which can only be rigorously defined within the framework of
distributions. Because of this term, the fluid limit equation
(see~\eqref{eq:soflu1}) is the integrated version of a partial differential equation rather than
an ordinary differential equation as it is the rule in the previously
studied queueing systems. Thus, the famous Gronwall's Lemma is of no
use here. Fortunately, the partial differential equation which pops
up, known as transport equation, is simple enough to have a
closed form solution -- see Theorem \ref{thm:sotran}. Thanks to that,
we can then proceed as usual to show the strong convergence of the
renormalized process to the fluid limit.

This paper is organized as follows. After some preliminaries, we define
and solve the integrated transport equation in the space of tempered
distributions. In Section  \ref{sec:ProTra}, we establish that the
above described process $\proc{\nu}$, is a weak Feller Markov process
and give its infinitesimal generator. In Section  \ref{sec:LGN}, we
prove the fluid limit theorem. The last section is devoted to
applications  to the EDF driven queue with deterministic initial time 
credits, and to a pure delay system.

\section{Preliminaries} 
\label{sec:PRELIM}
We denote by $\D_b$, respectively $\C_0$ and $\C_b$ the set of
real-valued functions defined on $\R$ which are bounded, right-continuous
with left-limits (rcll for short), respectively continuous vanishing
at infinity and bounded continuous. The space $\D_b$ is equipped with
the Skorokhod topology and   $\C_0$ and $\C_b$ with the topology of
the uniform convergence. The space of bounded differentiable
functions from $\R$ to itself is denoted by $\C^1_b$ and for $\phi\in
\C_b^1$, $\parallel \phi
\parallel_{\infty}:=\underset{x \in \R}{\sup}(\left|\phi(x)\right|+|\phi^\prime(x)|)$. 
For all $f \in \D_b$ and all $x \in \R$, we denote by $\tau_x f$, the function  
$\tau_xf(.):=f(.-x).$ 

The Schwartz space, denoted by $\S$, is the space of infinitely
differentiable functions, equipped with the topology defined by the
 semi-norms:
$$\mid \phi \mid_{a,b} :=\underset{x \in \R}{\sup} \mid x^a
\frac{d^b}{dx^b} \phi(x) \mid,\,a \in \N,\,b\in\N.$$ 
Its topological dual,  the space of tempered
distributions, is denoted by $\S^\prime$, and the duality product is classically denoted $\cro{\mu,\phi}$. The Fourier transform on $\S$ is defined by 
$\widehat \phi(\xi):=(2\pi)^{-1/2}\int_{\R}e^{-i\xi
  x}\phi(x)\,dx$ and the Fourier transform is defined on $\S^\prime$ by
the duality relation $\cro{\widehat{\mu}, \phi}=\cro{\mu,\, \widehat{\phi}}$. 

The set of finite positive measures on $\R$ is denoted by $\M_f^+$ 
and $\M_p$ is the set of finite counting measures on 
$\R$. The space $\M_f^+$ is embedded with the {weak} topology,  
$\sigma(\M_f^+,\C_b)$, 
for which $\M_f^+$ is Polish (we write $\cro{\mu,f}=\int f\,d\mu$ for $\mu \in \M_f^+$ and $f \in \D_b$).  We also denote for all $x\in \R$ and all 
$\nu \in \M_f^+$, $\tau_x\nu$ the measure satisfying for all 
Borel set $B$,  
$\tau_x\nu(B):=\nu\left(B-x\right).$  
Let $\C_0(\M_f^+,\R)$, be the set of continuous functions from  $\M_f^+$
to $ \R$, vanishing at  infinity, endpwed with the topology of the sup norm.
Let $0<T < \infty$, for $E$ a Polish space, we denote $\C\left([0,T],E\right)$, respectively $\D\left([0,T],E\right)$, the Polish space 
(for its usual strong topology) of continuous, respectively rcll,  
functions from $[0,T]$ to $E$.

\section{The integrated transport equation}
\label{subsec:TRANS}
The
\textit{transport equation} on 
$\mathcal \C^1(\R \times \R_+,\R)$ with unknown  $u(x,t)$ is defined as:
\begin{equation}
  \label{eq:tranport_usuel}
  \tag*{(E)}
  \begin{split}
  \partial_t u&=-b\partial_x u+f \text{ in } \R \times (0,\infty),\\
  u&=h \text{ at } \R \times \{t=0\},      
  \end{split}
\end{equation}
where $b$ is a real number, $f$ is a function of $\C^1\left(\R \times
  \R+,\R\right)$, and $h \in \C^1\left(\R,\R\right)$. It is well known
(see~\cite{MR99e:35001}) that \ref{eq:tranport_usuel} admits a unique
solution given for all $x,t$ by:
\begin{equation}
 u(x,t)=h(x-tb)+\int_0^t f(x+(s-t)b,s)ds.
\end{equation}
Let us define the following extension of the transport equation:
\begin{definition}
\textit{Let $T>0$, $K \in \S'$, $\proc{g} \in \D\left([0,T],\S'\right)$ such 
that $g_0\equiv 0$ and $b$ be a real number. 
The process $\proc{\eta}$ satisfies the integrated transport equation
E(K,g,b) on  $\D\left([0,T],\S'\right)$ if for all $\phi\in\S'$,
and for all $t\in [0,T]$: 
\begin{equation}
  \label{eq:Eg}
\tag*{E($K,g,b$)}
  \cro{\eta_t,\phi}=\cro{K,\phi}-b\int_0^t\cro{\eta_s,\phi^\prime}\,ds+\cro{g_t,\phi}.
\end{equation}
}
\end{definition}
\begin{theorem}
\label{thm:sotran}
The integrated transport equation {(\ref{eq:Eg})}  admits a unique
  solution $\proc{L}$ in $\D\left([0,T],\S'\right)$, satisfying for all $\phi
  \in \S$ and for all $t \in [0,T]$:
\begin{equation}
\label{solution0}
\cro{L_t,\phi}=\cro{K,\tau_{bt}\,\phi}+\cro{g_t,\phi}-b\int_0^t\cro{g_s,\tau_{b(t-s)}\,\phi^\prime}\,ds.
\end{equation}
\end{theorem}
\begin{proof}
Let $L$ and $M$ be two solutions of (\ref{eq:Eg}) and let $N=L-M$.
For all $t \in [0,T]$, it follows from \ref{eq:Eg} that for all $\phi \in \S$:
$$\frac{d}{dt}\cro{\widehat{N_t},\phi}
=-b\cro{N_t,\widehat \phi\ ^\prime}.$$
Denoting for all $\xi \in \R$, $\psi(\xi):=-i\xi$, this can be rewritten:
$$\frac{d}{dt}\cro{\widehat{N_t},\phi}=-b\cro{N_t,\widehat
  {\psi\phi}}=-b\cro{\widehat{N_t},\psi\phi}=-b\cro{\psi\widehat{N_t},\phi}.$$
Solving the latter differential equation yields  for all 
$\phi \in \S$ and all $t \in [0,T]$: 
$$\cro{\widehat{N_t},\phi}=\cro{\widehat{N_0}e^{b\psi
    t},\phi}=\cro{\widehat{N_0},e^{b\psi
      t}\phi}=\cro{N_0,\widehat{e^{b\psi t} \phi}}=0,$$
hence for all $t \in [0,T]$, $N_t\equiv  0$. Therefore, there is at most one solution to (\ref{eq:Eg}).

The process $\proc{L}$ defined by (\ref{solution0}) belongs to $\D\left([0,T],\S'\right)$ and for $\phi \in \S$, we have: 
\begin{multline*}
b\int_0^t\cro{L_s,\phi^\prime}\,ds\\=b\int_0^t\cro{K,\tau_{bs}
  \phi^\prime}\,ds+b\int_0^t\cro{g_s,\phi^\prime}\,ds-b^2\int_0^t\int_0^s\cro{g_r,\tau_{b(s-r)}\phi^\prime}\,dr\,ds.
\end{multline*} 
Since $\partial_t\cro{\zeta,\tau_{bt}\phi}=-b\cro{\zeta,\tau_{bt}\phi^\prime}$, we get:
\begin{multline*}
b\int_0^t\cro{L_s,\phi^\prime}\,ds=
-\int_0^t\frac{d}{ds}\left(\cro{K,\tau_{bs}\phi}\right)\,ds+b\int_0^t\cro{g_s,\phi^\prime}\,ds\\\shoveright{+b\int_0^t\int_r^t\frac{d}{ds}\left(\cro{g_r,\tau_{b(s-r)}\phi^\prime}\right)\,ds\,dr}\\
=-\cro{K,\tau_{bt}\phi}+\cro{K,\phi}+b\int_0^t\cro{g_s,\tau_{b(t-s)}\phi^\prime}\,ds
=-\cro{L_t,\phi}+\cro{K,\phi}+\cro{g_t,\phi},
\end{multline*}
The process $\proc{L}$ thus satisfies (\ref{eq:Eg}). 
\end{proof}

\section{The profile process}
\label{sec:ProTra}
Following Barrer's notation~\cite{MR19:779g}, we throughout this paper 
consider a
queueing system with impatient customers M/M/1/1+GI-EDF:  
\begin{itemize}
\item customers arrive at times $\suitei{T_i}$. The process defined for
  all $t$ by $$N_t:=\sum_{i\in\N^*}\car_{\left\{T_i\le t\right\}}$$ is 
  a Poisson process  of intensity $\lambda>0$,   
\item a first sequence of marks $\suitei{\sigma_i}$, the {sequence
    of service durations} requested by the customers, is i.i.d. with
  the distribution of $\sigma$ which an exponentially distributed with
  parameter $\mu > 0$ random variable, 
\item the customers are impatient: {i.e.}, the $i$-th customer
  leaves the system, and is lost forever, when he doesn't reach the
  service booth before his specific {deadline}, $T_i+D_i$.  In
  other words, he is initially labelled with a random variable
  referred to as his {patience}, or {initial time credit},
  $D_i$. The marks $\suitei{D_i}$ are independent and identically
  distributed with the distribution of $D$, an almost-surely
  non-negative and integrable random variable.  The time credits of
  the customers decrease continuously with time, at velocity one (in
  time units). Provided that the $i$-th customer entered the system
  before $t$ ($T_i \le t$), but did not reach the service booth before
  $t$, we denote $D_i(t)$ the {residual time credit at $t$} of
  this customer, {i.e.}, the residual time before his possible
  elimination.  Therefore:
$$D_i(t)=D_i-(t-T_i),$$
  and $D_i(t)\le 0$ means that the $i$-th customer has been
  {lost}, reaching his patience before $t$ before entering the service, 
\item there is $1$ non idling server and a buffer of infinite
  capacity, 
\item the service discipline is EDF ({i.e., Earliest Deadline
    First}): when completing a service, the server deals with the 
  customer whose residual time credit is the smallest among all the
  customers in the buffer, if any. This service then proceeds until
  completion, without any interruption.  
\end{itemize} 
Let us finally define the following {performance processes}: 
\[\left.
 \begin{array}{ll}
X_{t}&:= \text{ Number of customers in the system (buffer + service
  booth) at $t$ },\\
Q_{t}&:= (X_{t}-S)^+= \text{ Number of customers in the buffer at $t$ },\\
\S_{t}&:= \text{ Number of customers served up to time $t$ },\\
P_{t}&:= \text{ Number of customers lost up to time $t$ },
\end{array}\right.\]
 At time $t$, provided that the buffer is non-empty, denote for $i=1,....,Q_t$, 
$R_i(t)$ the $i$-th residual time credit of a customer in the buffer at $t$,
ranked in the increasing order: $$R_1(t)<R_2(t)<....<R_{Q_t}(t).$$
Provided that at least one customer has been lost at $t$ ($P_t \neq 0$), 
for $i=1,.....,P_t$, denote $R_{-i}(t)$, the $i$-th residual time credit
among the customers lost up to $t$ in the decreasing order: 
$$R_{-P_t}(t)<R_{-P_t-1}(t)<....<R_{-1}(t).$$
The {time credit profile} of the system at $t$ is the following measure:  
$$\nu_t:=\sum^{Q_t}_{i=1}\delta_{R_i(t)}+\sum^{P_t}_{i=1}\delta_{R_{-i}(t)},$$
where $\delta_x$ is  the Dirac mass at $x$. 
Provided that the buffer is non-empty at $t$, we denote for all $i=1,...,Q_t$, $$t_i(\nu_t):=R_i(t),$$ the $i$-th
point of $\nu_t$ (in the increasing order) on the positive half-line.

The service discipline can be represented as follows~:   
when the server completes a service (say at time $s$), he first deals with the 
customer whose time credit is given at this time by:
$$t_1(\nu_{s})=R_1(s)>0,$$
provided that $Q_{s} \neq 0.$ 
The customer corresponding to the atom $t_1(\nu_{s})$, being 
chosen by the server, leaves the buffer: the corresponding atom
$\delta_{t_1(\nu_s)}$ is erased from the point measure $\nu_t$
for all $t \ge s$ (this customer won't ever reappear in the buffer, since
the service discipline is non-preemptive).

\begin{figure}[!ht]
\begin{pspicture}(0,-2)(12,6.5)
  \psline(0,4)(10,4)
\psline(0,2.5)(10,2.5)
\psline(0,1)(10,1)
\put(11,4){$t$}
\put(11,2.5){$t+h$}
\put(11,1){$t+2h$}
\put(2.9,4.5){$0$}
\put(2.9,3){$0$}
\put(2.9,1.5){$0$}
\psline[linewidth=0.1](3,3.8)(3,4.2)
\psline[linewidth=0.1](3,2.3)(3,2.7)
\psline[linewidth=0.1](3,0.8)(3,1.2)
\psline(4,3.8)(4,4.2)
\psline(5.5,3.8)(5.5,4.2)
\psline(8,3.8)(8,4.2)
\psline(9.5,3.8)(9.5,4.2)
\put(3.7,4.5){$\scriptstyle R_1(t)$}
\put(5.3,4.5){$\scriptstyle R_2(t)$}
\put(7.7,4.5){$\scriptstyle R_3(t)$}
\put(9.2,4.5){$\scriptstyle R_4(t)$}
\psline(4,2.3)(4,2.7)
\psline(6.5,2.3)(6.5,2.7)
\psline(8,2.3)(8,2.7)
\put(3.8,2){$\scriptstyle R_1(t+h)$}
\put(6.2,2){$\scriptstyle R_2(t+h)$}
\put(7.7,2){$\scriptstyle R_3(t+h)$}
\psline(2.5,0.8)(2.5,1.2)
\psline(5,0.8)(5,1.2)
\psline(6.5,0.8)(6.5,1.2)
\psline(4,0.8)(4,1.2)
\psline[linestyle=dotted](2.5,1)(5.5,4)
\psline[linestyle=dotted](5,1)(8,4)
\psline[linestyle=dotted](6.5,1)(9.5,4)
\put(3.5,0.2){$\scriptstyle R_1(t+2h)$}
\put(4.5,-0.2){$\scriptstyle R_2(t+2h)$}
\put(6,0.2){$\scriptstyle R_3(t+2h)$}
\psline{->}(4,6)(4,5.5)
\psline(4,6)(5,6)
\put(5.3,6){\small{Beginning of service between $t$ and $t+h$}}
\psline{->}(2.5,0)(2.5,0.5){\small{Lost customer}}
\psline{->}(2,-0.5)(2,0)
\put(0,-1){\small{Arrival between $t$ and $t+2h$}}
\end{pspicture}
\caption{Dynamics of the profile process.}
\label{fig:profil}
\end{figure}
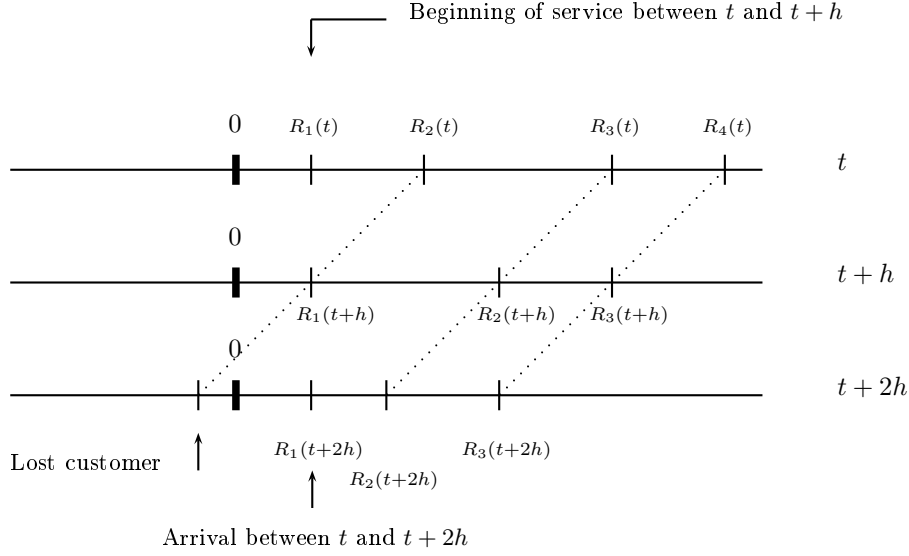

By {profile process  of the queue}, we mean, the process
$\proc{\nu}$ of the time credit profiles at $t$. 
This process is fully characterized by its initial value $\nu_0$, 
 the real numbers $\lambda > 0$ and $\mu > 0$ and  the non negative integrable random variable $D$. 
This process will consequently be referred to as the {profile 
  process associated to $\left(\nu_0,\lambda,\mu,D\right)$.} 
The dynamics of the profile process can be depicted as follows. The
atoms are translated continuously towards left at velocity 1,  
at the arrival time $T_i$, an atom is added to the measure
$\nu_{T_i}$ at
$D_i$ the initial time credit of the arriving customer, and at an end
of service $\tilde T_i$ , an atom disappear from the measure
$\nu_{\tilde T_i}$ at $t_1(\nu_{\tilde T_i}).$  
Figure \ref{fig:profil} shows a typical path of the profile process.
Note, that the buffer congestion and loss processes can be deduced from the profile process by writing for all 
$t\ge 0$: 
$$Q_t=\cro{\nu_t,\car_{\Rp}},\,
P_t=\cro{\nu_t,\car_{\R-}},$$
since the waiting, resp. already lost, customers at $t$ are those who have positive, resp. non positive, 
time credits at $t$. 

\section{Markov property}
\label{sec:markov-property}
Denote for all $t$, $A_t$ the remaining time before the next arrival
after $t$, and for all $t$ such that $X_t>0$, $
F_t$ the remaining time before the next end of service after $t$.    
For all $t,h>0$:
\[\nu_{t+h}=
 \left\{\begin{array}{ll}
   \tau_h\nu_t &\mbox{ if $A_t>h$ and $F_t>h$, or $Q_t=0$,}\\
   \tau_h\nu_t-\delta_{t_1(\nu_t)-h} &\mbox{ if $A_t>h$, $F_t<h$ and $Q_t>0$},\\
   \tau_h\nu_t+\delta_{d_k-(t+h-t_k)} &\mbox{ if $A_t<h$ and $F_t>h$ 
     or 
     $Q_t=0$}\\ &\mbox{ and the customer arrives at $t_k$ }\\&\mbox{ affected with
     the initial time credit $d_k$},
  \end{array}
  \right.\]
the more complex events (several arrivals, several ends of service,   
or arrivals and ends of service) between $t$ and $t+h$ being of probability $o(h)$. 
This dynamics shows in particular that $\proc{\nu} \in
\D\left([0,\infty),\M_f^+\right)$, since $\left(\cro{\nu_t,\phi}, t \ge
  0\right)$ belongs to $\D\left([0,\infty),\R\right)$ for all $\phi \in \C_b$.  
We finally
define the filtration:
$$\F_t:=\sigma\left(\nu_s(B),s\le t, B \in \B(\R)\right).$$ 
\begin{theorem}
\label{pro:Feller}
 The profile process $\proc{\nu}$ associated to
  $\left(\nu_0,\lambda,\mu,D\right)$ is a weak Feller process with
  respect to $\proc{\F},$ whose infinitesimal generator is given by:
\begin{multline}
\label{eq:generateur}
\A F(\nu)=\lim_{h
  \rightarrow 0}\frac{F(\tau_h\nu)-F(\nu)}{h}
-\Bigl(\lambda+\mu\car_{\left\{\nu_t(\Rp)>0\right\}}\Bigl)F(\nu)\\+\mu F\left(\nu-\delta_{t_1(\nu)}\right)\car_{\left\{\nu(\Rp)>0\right\}}
+\lambda\int F\left(\nu+\delta_d\right)\,d\P_D(d),
\end{multline}
for all $F$ in the {domain} of $\A$: 
$$\D(\A):= \C_0(\M_f,\R)\cap \left\{\lim_{h
  \rightarrow 0}\frac{F(\tau_h.)-F}{h} \mbox{
  exists} \right\}.$$
\end{theorem}
\begin{proof}
  For all $t,h\ge 0$ and all bounded measurable function
  $F:\M_f^+\rightarrow \R$:
  \begin{multline}
    \label{eq:transition}
    \begin{aligned}
      \esp{F(\nu_{t+h})|\F_t}=&\Bigl(1-\left(\lambda+\mu\car_{\left\{\nu_t(\Rp)>0\right\}}\right)h\Bigl)F(\tau_h\nu_t)\\&+\mu
      hF\left(\tau_h\nu_t-\tau_h\delta_{t_1(\nu_t)}\right)\car_{\left\{\nu_t(\Rp)>0\right\}}\end{aligned}\\
\shoveright{+\lambda
      h\int F\left(\tau_h\nu_t+\tau_h\delta_d\right)\,d\P_{D}(d) +
      o(h)}\\=:T_hF(\nu_t).
  \end{multline}
Thus, according to~\cite{Daw93}, p.18, $\proc{\nu}$ is a weak homogeneous
Markov process, whose transition function is given by $\left(T_h, h   \ge 0\right)$. 
For $F \in \C_0\left(\M_f^+,\R\right)$, it is easily seen from \req{eq:transition}, that  $T_hF \in
\C_0\left(\M_f^+,\R\right)$ for all 
$h\ge 0$. Since $\M_f^+$ embedded with the weak topology is 
  locally compact separable, it routinely follows that $\proc{\nu}$ is a weak Feller process whose infinitesimal generator of $\nu$ is given by \eqref{eq:generateur}.
\end{proof}

\begin{corollary}
\label{thm:martingale}
For all $\phi \in \C^1_b$, the process defined for all $t \ge 0$ by: 
\begin{multline}
\label{eq:mart1}
  M_{\phi}(t)=\cro{\nu_t,\phi}-\cro{\nu_0,\phi}-\int_0^t
  \cro{\nu_s,\phi^\prime}\,ds\\+\mu\int_0^t\phi\left(t_1(\nu_s)\right)\car_{\left\{\nu_s(\Rp)>0\right\}}\,ds-\lambda t \esp{\phi(D)}
\end{multline}
is an rcll $\F_t$-martingale such that $M_{\phi}(t) \in L^2$ for all
$t \ge 0$.  Its increasing process is given for all $t\ge
0$ by: 
\begin{equation}
\label{eq:croc1}
  <\!M_{\phi}\!>_t=\mu\int_0^t\phi^2\left(t_1(\nu_s)\right)\car_{\left\{\nu_s(\Rp)>0\right\}}\,ds+\lambda
t\esp{\phi^2(D)}.
\end{equation}
\end{corollary}
\begin{proof}
Let $\phi \in \C^1_b$. 
Define the mapping $\Pi_{\phi} :$  
$\M_f^+ \mapsto \R$ for all $\nu$ by: 
$$\Pi_{\phi}(\nu):=\cro{\nu,\phi}.$$
Since \begin{equation}
\label{eq:dualite2}
\lim_{h \rightarrow
  0}\frac{1}{h}\Bigl(\Pi_{\phi}(\tau_{h}\nu)-\Pi_{\phi}(\nu)\Bigl) 
=-\cro{\nu,\phi^\prime},
\end{equation} we have for all $\nu \in \M_p$: 
  \begin{equation*}
\A \Pi_{\phi}(\nu)=
-\cro{\nu,\phi^\prime} -\mu \phi\left(t_1(\nu)\right)\car_{\left\{\nu(\Rp)>0\right\}}
+\lambda\esp{\phi(D)}.
\end{equation*}
Furhtermore,
\begin{equation*}
\lim_{h \rightarrow
  0}\frac{1}{h}\Bigl(\Pi_{\phi}(\tau_{h}\nu)^2-\Pi_{\phi}(\nu)^2\Bigl)=-2\cro{\nu,\phi}\cro{\nu,\phi^\prime}, 
\end{equation*}
and hence for all $\nu \in M_p$:
\begin{equation*}
\A \Pi^2_{\phi}(\nu)
=2\cro{\nu,\phi}\A \Pi_{\phi}(\nu)+\mu\phi^2\left(t_1(\nu)\right)\car_{\left\{\nu(\Rp)>0\right\}}+\lambda\esp{\phi^2(D)}.
\end{equation*} 
From Dynkin's lemma~\cite{EthierKurtz,Dyn65}, it follows that
$M_\phi$ and the process defined for 
all $t$ by
\begin{equation*}
N_{\phi}(t)=\Pi^2_{\phi}(\nu_t)-\Pi^2_{\phi}(\nu_0)-\int_0^t\A
\Pi^2_{\phi}(\nu_s)\,ds
\end{equation*} are $\mathcal F_t$-local martingales. This entails that  
\begin{math}
  <\!M_{\phi}>_t=<\!\cro{\nu_.,\phi}\!>_t
\end{math}
for all $t \ge 0$, 
and thus, It\^o's integration by parts formula yields to:
\begin{multline*}
\cro{\nu_t,\phi}^2
=\cro{\nu_0,\phi}^2+2\int_0^t\cro{\nu_s,\phi}\,dM_{\phi}(t)\\+2\int_0^t\cro{\nu_s,\phi}
\A\Pi_{\phi}(\nu_s)\,ds+<\!\cro{\nu_.,\phi}\!>_t.
\end{multline*}
Hence, 
for all $t\ge 0$:
\begin{multline*}
2\int_0^t\cro{\nu_s,\phi}\,dM_{\phi}(t)+<\!\cro{\nu_.,\phi}\!>_t\\=N_{\phi}(t)+\mu\int_0^t\phi^2\left(t_1(\nu_s)\right)\car_{\left\{\nu_s(\Rp)>0\right\}}\,ds+\lambda
t\esp{\phi^2(D)}, 
\end{multline*}
and by identifying the finite variation processes, we obtain (\ref{eq:croc1}).
\end{proof}

\section{Fluid limit}
\label{sec:LGN}
For all $n \in \N^*$, denote by $\nu_0^{n}$, a measure on $\Rp$, 
and  
define $\proce{\nu_t^{n}}$, the profile 
process of the M/M/1/1+GI-EDF queue whose initial state is represented by the
profile $\nu_0^{n}$, whose arrival process $\proc{N^n}$ is Poisson of 
intensity 
$\lambda^n>$, where the customers request service durations 
 are exponentially distributed of mean expectation $(\mun)^{-1}$, and 
have initial
time credits i.i.d. with the distribution of $D^n$.  
The process $\proce{\nu_t^{n}}$ is in other words the profile process associated 
to $\left(\nu_0^n,\lambda^n,\mu^n,D^n\right)$.   
Also denote $\left(\mFnt, t \ge
  0\right)$, the associated filtration,
$$\ton:=\inf\left\{t \ge 0,
  \nu^n_t\left(\Rp\right)=0\right\},$$
the first time when the buffer is empty, and 
$$\omn:=\inf\left\{t \ge 0, t_1(\nu^n_t)=0\right\},$$
the first time of loss of the system. 
We also define as previously the performance processes of the $n$-th system: $\proc{X^n}$, $\proc{Q^n}$, given for all $t$ 
by $Q^n_t=\cro{\nu^n_t,\car_{\Rp}}$, $\proc{\S^n}$, $\proc{P^n}$, given by $P^n_t=\cro{\nu^n_t,\car_{R-}}$.\\ 
According to \thmref{thm:martingale}, 
for all $\phi \in \C^1_b$, the process defined for all $t \ge 0$ by: 
\begin{multline}
  \label{eq:martn1}
  M^n_{\phi}(t)=\cro{\nu^n_t,\phi}-\cro{\nu^n_0,\phi}-\int_0^t
  \cro{\nu^n_s,\phi^\prime}\,ds\\+\mu^n\int_0^t\phi\left(t_1(\nu^n_s)\right)\car_{\left\{\nu^n_s(\Rp)>0\right\}}\,ds-\lambda^n t \esp{\phi(D^n)}
\end{multline}
is an rcll $\mFnt$-martingale, such that $M^n_{\phi}(t) \in L^2$ for all $t$, and whose  
increasing process is given for all $t\ge 0$ by:
\begin{equation}
  \label{eq:crocn1}
  <M^n_{\phi}>_t=\mu^n\int_0^t\phi^2\left(t_1(\nu^n_s)\right)\car_{\left\{\nu^n_s(\Rp)>0\right\}}\,ds+\lambda^n
t\esp{\phi^2(D^n)}.
\end{equation}
We normalize the process $\proce{\nu^n_t}$ in time, space and weight the following way: for all Borel set $B$ and for all 
$t$, define  
$$\bar \nu^n_t(B)=\frac{\nu^n_{nt}(n B)}{n},$$
where $$n B:=\left\{nx, x \in B\right\}.$$ 
The first positive atom of $\bar \nu^n_t$ is
therefore given by:
$$t_1(\bar \nu^n_t)=\frac{t_1\left(\nu^n_{nt}\right)}{n}.$$
We also denote $\left(\mGnt, t \ge 0\right) :=
\left(\mathcal F^n_{nt}, t \ge 0\right)$, the associated filtration, 
$$\ovton:=\inf\left\{t \ge 0, \bar \nu^n_t
  (\Rp)=0\right\}=\frac{1}{n}\ton,$$
$$\ovomn:=\inf\left\{t\ge 0,t_1(\bar \nu^n_t)=0\right\}=\frac{1}{n}\omn$$ 
and normalize the arrival process as well as the performance processes of the $n$-th system the corresponding way, {i.e.}, for all $t \ge 0$, 
\begin{equation*}
\bar N^n_t:=\frac{N^n_{nt}}{t},\,\bar X^n_t:=\frac{X^n_{nt}}{t},\,\bar Q^n_t:=\frac{Q^n_{nt}}{t},\, 
\bar P^n_t:=\frac{P^n_{nt}}{t}. 
\end{equation*} 
For all $t\ge 0$, $\bar Q^n_t$ and $\bar P^n_t$ can thus be recovered by:
\begin{equation}
\label{eq:normQ}
\bar Q^n_t=\cro{\ovnun_t,\car_{\Rp}},
\end{equation}
\begin{equation}
\label{eq:normP}
\bar P^n_t=\cro{\ovnun_t,\car_{\R-}}.
\end{equation}
Let $\phi \in \C^1_b$ and $\psi^n(.)=\phi\left(./n\right)/n.$  
As easily seen from (\ref{eq:martn1}) and (\ref{eq:crocn1}),   
the process defined for all $t$ by
\begin{multline}
\label{eq:marnn1}
\barnp M(t):=M^n_{\psi^n}(nt)=\cro{\ovnun _t,\phi}-\cro{\ovnun _0,\phi}-\int_0^t \cro{\ovnun _s,\phi '}ds\\
+\mun \int_0^t \phi\left(t_1(\ovnun_s)\right)\car_{\{\ovnun_s(\Rp) > 0\}} \, ds - \lan  t\esp{\phi\left(\frac{D\n }{n}\right)}
\end{multline} 
is a $\mGnt$-martingale of $\D\left([0,\infty),\R\right)$, such that $\barnp M(t) \in L^2$ for all $t$.  
Its increasing process is given for all $t$ by:
\begin{equation*}
<\barnp M>_t\,=
\frac{\mun}{n}\int_{0}^{t}\phi^2\left(t_1(\ovnun_s)\right)\car_{\{\ovnun_s(\Rp) > 0\}}  \,ds+\frac{\lan}{n}t\esp{\phi^2\left(\frac{D\n}{n}\right)}.
\end{equation*}
We now define the set of hypothesis under which we will prove a law of
large numbers for the sequence of processes $\nu^n$.
\begin{hypo}
\label{hypo:H1}
{
\begin{itemize}
\item There exists two real numbers $\mu > 0$ and $\lambda > \mu$ such that:
\begin{equation}
      \label{eq:H12}
      \lan \tend \lambda,
    \end{equation}
\begin{equation}
      \label{eq:H13}
      \mun \tend \mu.
    \end{equation}
\item For all $\varepsilon>0$, there exists $M_{\varepsilon}>0$ such that for 
all $n\in \N^*$, 
\begin{equation}
\label{eq:H15bis}
\pr{\cro{\nu_0^{n},1} >nM_{\varepsilon}}\le \varepsilon.
\end{equation}
\item There exists a measure $\bar \nu_0^{*}$ of $\M_f^+$ such that for all $f \in \D_b$:
\begin{equation*}
      \suite{\cro{\ovnun_0,f}} \overset{\mathcal P}{\longrightarrow} \cro{\bar \nu_0^{*},f}.
    \end{equation*} 
\item There exist an integrable and almost surely non-negative r.v.
  $\bar D$ such that:
\begin{equation*}
      \frac{D^n}{n} \overset{\mathcal D}{\tend} \bar D,
    \end{equation*}
    \begin{equation*}
      \esp{\frac{D^n}{n}} \tend \esp{\bar D}.
    \end{equation*}
\end{itemize}
}
\end{hypo}
\begin{proposition}
Assume that Hypothesis \ref{hypo:H1} holds. 
Then, $\suite{\proce{\bar \nu^n_t}}$ is tight in $\D\left([0,\infty),\M_f^+\right)$.  
\end{proposition}

\begin{proof} 
According to Jakubowski's criterion~\cite{Daw93}, it suffices to show that:
\begin{enumerate}
\item For all $\phi \in \C^1_b$, the sequence $\suite{\proce{\cro{\ovnun_t,\phi}}}$ is tight in $\D\left([0,\infty),\R\right)$,
\item For all $T>0$ and $0<\eta<1$, there exists a compact subset $\mathbf K_{T,\eta}$ of $\M_f^+$ such that $$\underset{n \rightarrow \infty}{\lim \inf} \pr{\ovnun_t \in \mathbf K_{T,\eta} \forall t \in [0,T]}\ge 1-\eta.$$
\end{enumerate}
In order to prove the first condition, let us fix $\phi \in \C^1_b$ and $T>0$. 
Remarking that for all $s \ge 0$, $$\cro{\ovnun_s,1} \le \bar N^n_s + \cro{\ovnun_0,1},$$ 
Equation (\ref{eq:marnn1}) yields for all $u<v\le T$:
\begin{multline}
\label{eq:andreev0}
\left|\cro{\ovnun_v,\phi}-\cro{\ovnun_u,\phi}\right|
\le \int_u^v\left|\cro{\ovnun_s,\phi^\prime}\right|\,ds+\mu^n\int_u^v\left|\phi(t_1(\ovnun_s))\right|\,ds\\
\shoveright{+\lambda^n\esp{\phi\left(\frac{D^n}{n}\right)}|v-u|+\left|\barnp
  M(v)-\barnp M(u)\right|}\\
\shoveleft{\le |v-u|\parallel \phi^\prime
  \parallel_{\infty}\bar N^n_{T}+|v-u|\parallel \phi^\prime
  \parallel_{\infty}\cro{\ovnun_0,1}}\\
+|v-u|\parallel\phi\parallel_{\infty}\left(\left|\mun-\mu\right|+\left|\lan-\lambda\right|\right)
+|v-u|\parallel\phi\parallel_{\infty}\left(\lambda+\mu\right)\\+\left|\barnp
  M(v)-\barnp M(u)\right|.
\end{multline}
Let $\varepsilon>0$ and $\eta>0$. First, 
let
$$\delta_{1}:=\frac{\varepsilon\eta}{30\parallel\phi^\prime\parallel_{\infty}\lambda
T}.$$
From Markov's inequality:
\begin{multline*}
\pr{\underset{u,v<T\,,\,\mid v-u\mid\le\delta_{1}}{\sup}|v-u|\parallel \phi^\prime
  \parallel_{\infty}\bar N^n_{T} \ge \frac{\eta}{5}}\\
\le \frac{5\delta_{1}\parallel\phi^\prime\parallel_{\infty}}{\eta}\esp{\bar
  N^n_{T}}\\
\le \frac{\lambda T
  5\delta_{1}\parallel\phi^\prime\parallel_{\infty}}{\eta}+\frac{\left|\lambda-\lambda^n\right|T 5\delta_{1}\parallel\phi^\prime\parallel_{\infty}}{\eta}\\=\frac{\varepsilon}{6}+\frac{\left|\lambda-\lambda^n\right|T 5\delta_{1}\parallel\phi^\prime\parallel_{\infty}}{\eta}, 
\end{multline*}
and thus with (\ref{eq:H12}), there exists $N_{1}>0$ such that for all
$n \ge N_{1}$, 
\begin{equation*}
\pr{\underset{u,v<T\,,\,\mid v-u\mid\le\delta_{1}}{\sup}|v-u|\parallel \phi^\prime
  \parallel_{\infty}\bar N^n_{T} \ge \frac{\eta}{5}}\le \frac{\varepsilon}{6}+\frac{\varepsilon}{6}=\frac{\varepsilon}{3}.
\end{equation*} 
Let now:
$$\delta_2:=\frac{\eta}{5M_{\varepsilon/3}\parallel\phi^\prime\parallel_{\infty}}.$$
According to (\ref{eq:H15bis}), for all $n \in \N^*$, 
\begin{multline*}
\pr{\underset{u,v<T\,,\,\mid v-u\mid\le\delta_{2}}{\sup}|v-u|\parallel \phi^\prime
  \parallel_{\infty}\cro{\ovnun_0,1} \ge \frac{\eta}{5}}\\
\le \pr{\cro{\ovnun_0,1}\ge \frac{\eta}{5\delta_2\parallel \phi^\prime
  \parallel_{\infty}}}
=\pr{\cro{\ovnun_0,1}\ge M_{\varepsilon/3}}
\le \frac{\varepsilon}{3}.
\end{multline*}
According to assumptions (\ref{eq:H12}) and (\ref{eq:H13}), there exists $N_2$,
such that for all $n \ge N_2$, for all $u,v \le T$,
\begin{equation*}
|v-u|\parallel\phi\parallel_{\infty}\left(\left|\mun-\mu\right|+\left|\lan-\lambda\right|\right)
< \frac{\eta}{5},
\end{equation*}
and letting 
$$\delta_3:=\frac{\eta}{6\parallel\phi\parallel_{\infty}(\lambda+\mu)},$$
\begin{equation*}
  \underset{u,v<T\,,\,\mid
    v-u\mid\le\delta_{3}}{\sup}|v-u|\parallel\phi\parallel_{\infty}\left(\lambda+\mu\right)\le \delta_3\parallel\phi\parallel_{\infty}\left(\lambda+\mu\right) < \frac{\eta}{5}.
\end{equation*}
Now, let $\xi>0$. Apply successively Markov's and Doob's inequalities:
\begin{multline*}
\pr{\underset{t \le T}{\sup}\left|\barnp M(t)\right| \ge \xi}
\le\frac{4}{\xi^2}\esp{<\barnp
      M>_{T}}\\
=\frac{4}{\xi^2}\esp{\frac{\mun}{n}\int_{0}^{T}\phi^2\left(\mathcal
    R(\ovnun_s)\right)\car_{\left\{\ovnun_s(\Rp) >0\right\}}\,ds+\frac{\lan}{n}T\esp{\phi^2\left(\frac{D\n}{n}\right)}}\\
\le
\frac{4}{\xi^2}\left(\frac{\mu^n}{n}+\frac{\lambda^n}{n}\right)\mid\mid\phi^2\mid\mid_{\infty}T
\tend 0.
\end{multline*}
For all $n\in\N^*$, $\procp{\barnp M}$ being a rcll process on
$[0,T]$, one can apply the standard convergence criterion~\cite{Rob}, from which it follows that 
$\suite{\procp{\barnp M}}$ converges in distribution to the null process. This sequence is in particular tight in $\D\left([0,T],\R\right)$:
there exists $\delta_4>0$ and $N_3>0$ such that for all $n \ge N_3$: 
\begin{equation*}
  \pr{\underset{u,v \le T\,,\,|v-u|\le \delta_4}{\sup}\left|\barnp
  M(v)-\barnp M(u)\right| \ge \frac{\eta}{5}}\le \frac{\varepsilon}{3}.
\end{equation*}
Finally, in view of the previous inequalities and (\ref{eq:andreev0}), 
there exists $\delta >0$ and 
$N \in \N$ such that for all 
$n \ge N$:
\begin{equation}
\label{eq:tightness2}
\pr{\underset{u,v \le T\,,\,\left|v-u\right|\le
    \delta}{\sup}\left|\cro{\ovnun_v,\phi}-\cro{\ovnun_u,\phi}\right|\ge \eta}\\
\le\varepsilon.
\end{equation}
On the other hand, let $$\alpha_{\varepsilon}:=M_{\varepsilon}\parallel\phi\parallel_{\infty}.$$ 
Assumption (\ref{eq:H15bis}) implies that 
for all $n \in \N^*$, 
\begin{equation}
  \label{eq:tightness1}
  \pr{\left|\cro{\ovnun_0,\phi}\right|>\alpha_{\varepsilon}}
\le \mathbf P\biggl[\parallel\phi\parallel_{\infty}\cro{\ovnun_0,1}>M_{\varepsilon}\parallel\phi\parallel_{\infty}\biggl]\le \varepsilon.
\end{equation}
With (\ref{eq:tightness1}) and (\ref{eq:tightness2}) we can apply the
standard tightness criterion of real valued processes (see for instance ~\cite{Rob}): for all $T>0$,
$\suite{\proce{\cro{\bar \nu^n_t,\phi}}}$ is tight in $\D\left([0,T],\R\right)$: 
it is tight in $\D\left([0,\infty),\R\right)$.  

 We now
prove the second tightness condition (compact containment).     
Let us first apply ~\cite{Gromoll}, Lemma A.2.: under hypothesis 
\ref{hypo:H1}, we have the following weak law of large numbers: 
$$\suite{\proce{\frac{1}{n}\sum_{i=1}^{\bar N^n_t}\phi(\bar
    D^n_i)}}\Longrightarrow \proce{\lambda t\esp{\phi(\bar D)}} \mbox{ in }\D\left([0,T],\R\right),$$
for any $\phi \in \C_b$. In particular, this yields for any $0<l\le
T$:
\begin{equation}
\label{eq:tight1}
\pr{\underset{t \in [0,T-l]}{\sup}\frac{1}{n}\sum_{N^n_{nt}+1}^{N^n_{n(t+l)}}\phi(\bar D^n_i) > 2\lambda
    l\esp{\phi(\bar D)}}\tend 0.
\end{equation} 
Taking $l=T$ and $\phi=1$ in the last expression
yields:
\begin{equation*}
  \pr{\bar N^n_T >2 \lambda T}\tend 0.
\end{equation*}
Denote $I(.)$, the identity on $\R$. Taking $l=T$ and $\phi=I$ in (\ref{eq:tight1}) also leads to: 
\begin{equation*}
  \pr{\frac{1}{n}\sum_{i=1}^{N_{nT}^n}\bar D^n_i>2\lambda T\esp{\bar D}}\tend 0.
\end{equation*}
Let $$M_T=\max\biggl\{2\lambda T+\cro{\bar \nu_0^{* },1},2\lambda T\esp{\bar D}+\cro{\bar\nu_0^{* },I}\biggl\}+1.$$ 
We have:
\begin{multline}
\label{eq:tight4}
\pr{\underset{t \in [0,T]}{\sup}\max\bigl\{\cro{\bar \nu^n_t,\car_{\R+}},\cro{\bar \nu^n_t,I\car_{\R+}}\bigl\}>M_T}\\
\shoveleft{\le\pr{\cro{\bar \nu_0^{n},1}+\bar
N^n_T>2\lambda T+\cro{\bar\nu_0^{* },1}+1}}\\\shoveright{+\pr{\cro{\bar \nu_0^{n},I}+\frac{1}{n}\sum_{i=1}^{N^n_{nT}}\bar D^n_i>2\lambda T\esp{\bar D}+\cro{\bar\nu_0^{* },I}+1}}\\
\shoveleft{\le\pr{\cro{\bar \nu_0^{n},1}>\cro{\bar\nu_0^{* },1}+1}+\pr{\bar
N^n_T>2\lambda T}}\\\shoveright{+\pr{\cro{\bar \nu_0^{n},I}>\cro{\bar\nu_0^{* },I}+1}+
\pr{\frac{1}{n}\sum_{i=1}^{N_{nT}^n}\bar D^n_i>2\lambda T\esp{\bar D}}}\\\tend 0.
\end{multline}
\textit{ }\\
Let us now define, for all $T>0$, $0<\eta<1$, the set $$\mathcal
K_{T,\eta}:=\biggl\{\zeta \in
  \M_f^+\,;\,\max\bigl\{\cro{\zeta,\car_{\R+}},\cro{\zeta,I\car_{\R+}}\bigl\}\le
  M_T\,,\,\cro{\zeta,\car_{\left(-\infty,-T\right]}}=0\biggl\}.$$ 
Since $\cro{\zeta,I\car_{\R+}}\le M_T$, this  implies that for all $y>0$,
$\zeta\left([y,\infty)\right)\le M_T/y,$
and thus 
$$\underset{y \rightarrow \infty}{\lim} \underset{\zeta \in \mathcal
  K_{T,\eta}}{\sup}\zeta\left([y,\infty)\right)=0\,\,\,,\,\,\,\underset{y \rightarrow -\infty}{\lim} \underset{\zeta \in \mathcal
  K_{T,\eta}}{\sup}\zeta\left(\left(-\infty,y\right])\right)=0,$$
which implies that $\mathcal K_{T\eta} \subset \M_f^+$ is relatively
compact~\cite{kallenberg83}. Now, since up to time $T$ no lost 
customer can have a residual time credit less than $-T$, 
$$\underset{t \le
  T}{\sup}\cro{\bar \nu^n_t,\car_{\left(-\infty,-T\right]}}=0.$$ This, together
with (\ref{eq:tight4}) implies that:
$$\underset{n \rightarrow \infty}{\liminf}\mathbf P\biggl[\bar \nu^n_t \in \mathcal
  K_{T,\eta}, \mbox{ for all }t \in [0,T]\biggl]\ge 1-\eta. $$ 
$\mathbf{K}_{T,\eta}$ being the closure of $\mathcal K_{T,\eta}$,
we found a compact subset $\mathbf {K}_{T,\eta}
\subset \M_f^+$ such that: 
\begin{equation*}
  \underset{n \rightarrow \infty}{\liminf}\mathbf P\biggl[\bar \nu^n_t \in \mathbf
  {K}_{T,\eta}, \mbox{ for all }t \in [0,T]\biggl] \ge 1-\eta.
\end{equation*}
\end{proof}
\begin{theorem}[Fluid limit theorem for M/M/1/1+GI-EDF queues]
\label{thm:oui}
Assume that Hypothesis \ref{hypo:H1} holds, that there exists $T>0$ such that:
\begin{equation}
\label{eq:TpsVidage}
\pr{\ovton \le T} \tend 0,
\end{equation}
and that there exists a 
deterministic element $\proce{\ovrt}$ of $\C\left([0,\infty),\R\right)$ such that:
\begin{equation}
\label{eq:ConvPre}
\proce{t_1(\bar \nu^n_t)} \Longrightarrow \proce{\ovrt} \mbox{ in $\D\left([0,T],\R\right)$}.
\end{equation}
 Then:
$$\suite{\proce{\bar \nu^n_t}} \Longrightarrow 
\proce{\bar \nu^*_t}\,\mbox{ in } \D\left([0,T],\M_f^+\right),$$ 
where $\proce{\bar \nu^*_t}$ is the deterministic element of $\C\left([0,T],\M_f^+\right)$ defined for all $f \in \D_b$ and all $t \in [0,T]$ by:
\begin{equation*}
\cro{\bar \nu_t^{* },f}=\cro{\bar \nu_0^{* },\tau_{t}\,f}-\mu\int_0^t \tau_{t-s}f\left(\bar
  r_{s} \right)\,ds
+\lambda \int_0^t \esp{\tau_{t-s}f\left(\bar D\right)}\,ds.
\end{equation*} 
\end{theorem}
\begin{proof} 
Let $\proc{\bar \chi}$ be a limit point of $\suite{\proce{\bar \nu^n_t}}$. 
On the one hand, assumption (\ref{eq:TpsVidage}) implies 
that  $\suite{\car_{\left\{\ovton \le T\right\}}}$ converges in distribution to $0$. 
On another hand, for all $\phi \in \C^1_b$ the mappings 
\[\Psi_1:\left\{
\begin{array}{ll}
\D\left([0,\infty),\R\right)&\hookrightarrow \D\left([0,\infty),\R\right)\\
\proc{X}&\mapsto \proce{\phi(X_t)}
\end{array}\right.\]
and 
\[\Psi_2:\left\{
\begin{array}{ll}
\displaystyle\D\left([0,\infty),\R\right)&\hookrightarrow \C\left([0,\infty),\R\right)\\
\displaystyle\proc{Y}&\mapsto \proce{\int_0^t Y_s\,ds}
\end{array}\right.\]
are continuous, as well as $\Psi:=\Psi_2 \circ \Psi_1,$ hence in view 
of (\ref{eq:ConvPre}), the
continuous mapping theorem entails that:
$$\suite{\proce{\int_0^t\mathcal \phi\left(t_1(\ovnun_s)\right)\,ds}}
\Longrightarrow \proce{\int_0^t\phi(\bar r_s )\,ds}\mbox{ in }\D\left([0,T],\R\right).
$$  
Consequently, for
all $t \in [0,T]$, all $\phi\in \C^1_b$:
\begin{equation}
\label{eq:soflu1}
\cro{\bar \chi_t,\phi}=\cro{\bar\nu_0^{* },\phi}-\int_0^t\cro{\bar
  \chi_s,\phi^\prime}\,ds-\mu\int_0^t \phi\left(\bar r_s \right)\,ds+\lambda t\esp{\phi(\bar D)}.
\end{equation} 
In particular, the latter is true for all $\phi \in \S$: 
(\ref{eq:soflu1}) is the integrated transport equation (E$(\bar\nu_0^{* }, g, 1)$), 
where $g$ is defined by: for all $\phi\in {\mathcal S}$,
$$\cro{g_t,\phi}:=-\mu\int_0^t \phi\left(\bar r_s \right)\,ds+\lambda t\esp{\phi(\bar D)}.$$
According to \thmref{thm:sotran}, the only solution of
(\ref{eq:soflu1}) is given for all $t \in[0,T]$ and all $\phi \in \S$ by:
\begin{multline*}
\cro{\bar \chi_t,\phi}=\cro{\bar\nu_0^{* },\tau_{t}\,\phi}+\cro{g_t,\phi}-\int_0^t\cro{g_s,\tau_{t-s}\,\phi^\prime}\,ds\\
\shoveleft{=\cro{\bar\nu_0^{* },\tau_{t}\,\phi}-\mu\int_0^t \phi\left(\bar
  r_{s} \right)\,ds+\lambda t\esp{\phi(\bar D)}}\\
+\mu\int_0^t\int_0^s \left(\tau_{t-s}\phi^\prime\left(\bar
  r_{u} \right)\right)\,du\,ds-\lambda\int_0^t s\esp{\tau_{(t-s)}\phi^\prime\left(\bar D\right)}\,ds\\
=\cro{\bar \nu^*_t,\phi}.
\end{multline*}
The limit point is therefore unique in $\D\left([0,T],\M_f^+\right)$, equal to
$\proce{\bar \nu^*_t}$, since $\S$ is a separating class of $\M_f^+$. 
\end{proof}

\section{Applications}
\label{sec:DET}

\subsection*{{M/M/1/1+D-EDF case}}
We hereafter apply Theorem \ref{thm:oui} to determine the fluid limit 
of the M/M/1/1+GI-EDF system in which the time credits of
the customers are deterministic. 
 We verify in particular that
assumptions (\ref{eq:TpsVidage}) and (\ref{eq:ConvPre}) are satisfied
in this case, and specify the form of the limit.

We therefore consider a sequence of M/M/1/1+D-EDF systems, for which we make the following assumptions: 
\begin{hypo}[Basic Assumptions for a M/M/1/1+D-EDF system]
\label{hypo:HDET}
{
\begin{itemize}
\item For all $n \in \N^*$, there are initially $n+1$ customers in the buffer, all of them with time 
credit $nd$, where $d>0$ (that is, the $n$ customers have their deadline at time $nd$), 
\item for all $n \in \N^*$, $\lan$ is the intensity of the Poisson process of arrivals, where $$\lan \tend \lambda>0,$$ 
\item for all $n \in \N^*$  the customers require service durations 
exponentially distributed,
  of parameter $\mun$, satisfying: $$\mun \tend \mu\,,\mbox{ where }(d)^{-1}<\mu<\lambda,$$ 
\item for all $n \in \N^*$, the initial time credit of any customer is 
deterministic, given by $d^n$, where $d^n/n 
  \tend d.$ 
\end{itemize}}
\end{hypo}
This queueing system is described by the profile process 
$\proce{\bar \nu_t^{n,\text{\tiny{D}}}}$, associated to $\left(n\delta_{nd},\lambda^n,\mu^n,d^n\right)$, which keeps track of all the residual time credits of
all the customers waiting in the buffer are already lost.
The notations are those of the preceding sections, with superscripts 
$^{\text{\tiny{D}}}$, for ``deterministic''.  


\begin{lemma}
\label{lemma:Vidage}
For all $x < \mu^{-1}$, 
$$\pr{\bar \tau_0^{n,\text{\tiny{D}}} \le x} \tend 0.$$
\end{lemma}

\begin{proof}
For the event
$$\left\{\bar \tau_0^{n,\text{\tiny{D}}} \le x\right\}=\left\{\tau_0^{n,\text{\tiny{D}}} \le nx\right\}$$
to occur, the $n+1$ customers initially present in the buffer must have all entered 
the service before time $nx$,  since they couldn't have been eliminated before $nx \le n/\mu < nd.$ Therefore the first $n$ customers among them must have completed their service before $nx$, or in other words:
\begin{equation*}
\pr{\bar \tau_0^{n,\text{\tiny{D}}} \le x} = \pr{\tau_0^{n,\text{\tiny{D}}} \le nx}\le \pr{S^n_{nx} \ge n},
\end{equation*} 
where $\proce{S^n_t}$ denotes a Poisson process of intensity $\mu^n$
(the server works without interruption at least until he has completed
the services of these customers). Hence, $S^n_{nx}$  has the same distribution as the sum of $n$ independent r.v $\left(P^n_i,i=1,...,n\right)$, Poisson distributed of parameter $\mu^nx$. Hence, denoting $\left(P_i,i=1,...,n\right)$, a family of $n$ independent r.v. Poisson distributed of parameter $\mu$,
\begin{multline*}
\pr{\bar \tau_0^{n,\text{\tiny{D}}} \le x}\le\pr{\frac{1}{n}\sum_{i=1}^n P^n_i \ge 1}\underset{\infty}{\sim} \pr{\frac{1}{n}\sum_{i=1}^n P_i \ge 1}\tend \car_{\left\{x\mu \ge 1\right\}}=0, 
\end{multline*}
according to the weak law of large numbers.
\end{proof}

\begin{lemma}
\label{lemma:convomega}
For all $\xi >0$, $$\pr{\bar \omega_0^{n,\text{\tiny{D}}} \le \bar \omega_0^{*,\text{\tiny{D}}}-\xi}
\tend 0,$$
where $$\bar \omega_0^{*,\text{\tiny{D}}}:=\frac{\rho d-\mu^{-1}}{\rho-1}.$$
\end{lemma}

\begin{proof}
We have:
\begin{multline}
\label{eq:bach1}
\pr{\bar \omega_0^{n,\text{\tiny{D}}} \le \bar
  \omega_0^{*,\text{\tiny{D}}}-\xi}
=\pr{\left\{\omega_0^{n,\text{\tiny{D}}} \le n\bar \omega_0^{*,\text{\tiny{D}}}-n\xi\right\}\cap\left\{n\bar \omega_0^{*,\text{\tiny{D}}}-n\xi\le
    \tau_0^{n,\text{\tiny{D}}}\right\}}\\+
\pr{\left\{\omega_0^{n,\text{\tiny{D}}} \le \tau_0^{n,\text{\tiny{D}}}\right\}\cap\left\{n\bar \omega_0^{*,\text{\tiny{D}}}-n\xi>
    \tau_0^{n,\text{\tiny{D}}}\right\}}\\
\shoveright{+\pr{\left\{\tau_0^{n,\text{\tiny{D}}}\le\omega_0^{n,\text{\tiny{D}}} \le n\bar \omega_0^{*,\text{\tiny{D}}}-n\xi\right\}\cap\left\{n\bar \omega_0^{*,\text{\tiny{D}}}-n\xi>
    \tau_0^{n,\text{\tiny{D}}}\right\}}}\\
\le \pr{\omega_0^{n,\text{\tiny{D}}} \le (n\bar \omega_0^{*,\text{\tiny{D}}}-n\xi)\wedge \tau_0^{n,\text{\tiny{D}}}}+\pr{\tau_0^{n,\text{\tiny{D}}}\le\omega_0^{n,\text{\tiny{D}}}\le n\bar \omega_0^{*,\text{\tiny{D}}}-n\xi}.
\end{multline}
 Let us denote for all $t \ge 0$:
\[\left\{\begin{array}{ll}
\mathcal N^n_t:=&\mbox{number of customers arrived up to
    time $t$,}\\
&\,\,\,\mbox{having deadline before $t$}\\
\mathcal S^n_t:=&\mbox{number of services completed up to
    $t$}.
\end{array}\right.\]
On the event $$\left\{\omega_0^{n,\text{\tiny{D}}} \le (n\bar
  \omega_0^{*,\text{\tiny{D}}}-n\xi)\wedge
  \tau_0^{n,\text{\tiny{D}}}\right\},$$ there is at least one loss and
no idle time before time $(n\bar
\omega_0^{*,\text{\tiny{D}}}-n\xi)\wedge\tau_0^{n,\text{\tiny{D}}}$. For
this to occur, since the service discipline amounts to FIFO, there
must be at the first time of loss, say $nt \ge nd$, the number of
customers initially in the system or arrived up to $nt$, of a priority higher to the priority of the customer who is
lost at $nt$ (there are
$n+1+\mathcal N^n_{nt}$ such customers) must be greater than the
number of services initiated up to $nt$ (i.e., $\mathcal
S^n_{nt}+1$). We can therefore write that:
\begin{multline}
\label{eq:spleen1}
\pr{\omega_0^{n,\text{\tiny{D}}} \le (n\bar \omega_0^{*,\text{\tiny{D}}}-n\xi)\wedge \tau_0^{n,\text{\tiny{D}}}}\\\le
\pr{\underset{d \le t \le (\bar \omega_0^{*,\text{\tiny{D}}}-\xi)\wedge\bar \tau_0^{n,\text{\tiny{D}}}}{\sup} 
  \left(\mathcal N^n_{nt}+n+1\right)-\left(\mathcal S^n_{nt}+1\right)\ge 0}\\
=
\pr{\underset{d \le t \le (\bar \omega_0^{*,\text{\tiny{D}}}-\xi)\wedge\bar \tau_0^{n,\text{\tiny{D}}}}{\sup} \mathcal
  N^n_{nt}-\mathcal S^n_{nt}\ge -n }.
\end{multline}
Since there has been no idle time in $\left[0,n(\bar \omega_0^{*,\text{\tiny{D}}}-\xi)\wedge\tau_0^{n,\text{\tiny{D}}}\right]$, for all $t$ in this interval, $\mathcal S^n_{nt}$ has the same distribution as $S^n_{nt},$  
where $\proc{S^n}$ is a Poisson process of intensity $\mun$. 
On the other hand, the process of arrivals $\proc{N^n}$ marked by the
initial time credits of the customers $\suitei{D^n_i}$ being a
two-dimensional Poisson process, it is easily checked that the process 
$\proce{\mathcal N^n_{t}-\lambda^n_t}$ is an $\mathcal
F^n_t$-martingale, where 
$$\lambda^n_t=\lan\left(t-n d\right)^+.$$
Thus the process defined for all $t$ by  
$$\mathcal M^n_t:=\mathcal N^n_{nt}-
S^n_{nt}-\left(\lambda^n_{nt}-\mu^nnt\right)$$
is a $\mathcal G^n_t$-martingale.  
Hence, with (\ref{eq:spleen1}): 
\begin{multline}
\label{eq:bach2}
\pr{\omega_0^{n,\text{\tiny{D}}} \le (n\bar \omega_0^{*,\text{\tiny{D}}}-n\xi)\wedge \tau_0^{n,\text{\tiny{D}}}}\\
\le \pr{\underset{d \le t \le \bar \omega_0^{*,\text{\tiny{D}}}-\xi}{\sup} \mathcal M^n_{t}\ge \underset{d \le t \le \bar \omega_0^{*,\text{\tiny{D}}}-\xi}{\inf}\left(\mu^nnt-\lambda^n_{nt}\right)-n}\\
=\pr{\underset{t \le \bar \omega_0^{*,\text{\tiny{D}}} -
    \xi}{\sup} \mathcal M^n_{t}\ge n\left\{\lambda^nd-1-\left(\lambda^n-\mu^n\right)\left(\bar \omega_0^{*,\text{\tiny{D}}}-\xi\right)\right\}}\\
\le \frac{4}{n^2\left\{\lambda^nd-1-\left(\lambda^n-\mu^n\right)\left(\bar \omega_0^{*,\text{\tiny{D}}}-\xi\right)\right\}^2}\esp{\cro{\mathcal
      M^n}_{\bar \omega_0^{*,\text{\tiny{D}}}-\xi}}\\
\underset{n \rightarrow \infty}{\sim}\frac{4}{\left(n(\lambda-\mu)\xi\right)^2}\left\{n\lambda^n\left(\bar \omega_0^{*,\text{\tiny{D}}} -
    \xi-d\right)+\mu^n
  n\left(\bar \omega_0^{*,\text{\tiny{D}}}-\xi\right)\right\}\\\tend 0,
\end{multline}
using successively Doob's inequality and the fact that $(\lambda-\mu)\bar \omega_0^{*,\text{\tiny{D}}}=\lambda d-1.$
Now, clearly 
\begin{multline}
\label{eq:Kat1}
\pr{\tau_0^{n,\text{\tiny{D}}}\le\omega_0^{n,\text{\tiny{D}}} \le
  n\bar \omega_0^{*,\text{\tiny{D}}}-n\xi}\\\le\pr{\frac{n}{2\mu}\le \tau_0^{n,\text{\tiny{D}}}\le\omega_0^{n,\text{\tiny{D}}} \le
  n\bar \omega_0^{*,\text{\tiny{D}}}-n\xi}+\pr{\tau_0^{n,\text{\tiny{D}}} \le \frac{n}{2\mu}}.
\end{multline}
on the event $$\left\{\frac{n}{2\mu}\le\tau_0^{n,\text{\tiny{D}}}\le\omega_0^{n,\text{\tiny{D}}} \le
  n\bar \omega_0^{*,\text{\tiny{D}}}-n\xi\right\},$$ there exists 
$t \in \left[\frac{1}{2\mu}, \bar \omega_0^{*,\text{\tiny{D}}}-\xi\right]$ (the first one) such that the
buffer is empty at $nt$, but there has been no loss before $nt$. For this event
to occur, there must be up to $nt$, the same number of customers
entered ({i.e.}, $N^n_{nt}$) as of services initiated ({i.e.}, $\mathcal S^n_{nt}+1$). 
Thus, remarking that the process defined for all $t$ by:
$$M^n_{t}:=N^n_{nt}-S^n_{nt}-(\lan nt - \mun nt)$$
is a $\mathcal G^n_t$-martingale, and that for all $t\le \tau_0^{n,\text{\tiny{D}}}$, 
$\mathcal S^n_{nt}$ equals $S^n_{nt}$ in distribution,  
\begin{multline*}
\pr{\frac{n}{2\mu}\le\tau_0^{n,\text{\tiny{D}}}\le\omega_0^{n,\text{\tiny{D}}} \le n\bar \omega_0^{*,\text{\tiny{D}}}-n\xi}\\ \le
\pr{\underset{1/(2\mu) \le t \le
    \bar \omega_0^{*,\text{\tiny{D}}}-\xi}{\sup}\mathcal S^n_{nt}-N^n_{nt} \ge -1}\\
\le\pr{\underset{1/(2\mu) \le t \le
    \bar \omega_0^{*,\text{\tiny{D}}}-\xi}{\sup}-M^n_{t} \ge
  -1+n(\lan-\mun)\frac{1}{2\mu}}\\
\le \frac{4}{\left(n(\lan-\mun)\frac{1}{2\mu}-1\right)^2}\esp{<M^n>_{\bar \omega_0^{*,\text{\tiny{D}}}-\xi}}\\
\underset{n \rightarrow \infty}{\sim} \frac{16}{\left(n(\rho-1)-2\right)^2}\left\{\left(\lambda^n
    +\mun\right) n(\bar \omega_0^{*,\text{\tiny{D}}}-\xi)\right\} 
\tend 0, 
\end{multline*}  
which, together with Lemma \ref{lemma:Vidage} applied to $x=1/(2\mu)$ and (\ref{eq:Kat1}) yields:
\begin{equation}
  \label{eq:bach3}
  \pr{\tau_0^{n,\text{\tiny{D}}}\le\omega_0^{n,\text{\tiny{D}}} \le n\bar \omega_0^{*,\text{\tiny{D}}}-n\xi}\tend 0.
\end{equation}
We conclude by substituting (\ref{eq:bach2}) and (\ref{eq:bach3}) in (\ref{eq:bach1}).
\end{proof}

\begin{proposition}
\label{pro:ConvPre}
Assume Hypothesis 1 and 2 holds, then 
$$\suite{\proce{t_1(\bar \nu_t^{n,\text{\tiny{D}}})}} \Longrightarrow 
\proce{\bar r_t^{\text{\tiny{D}}}} \mbox{ in } \D\left([0,\infty),\R\right),$$
where $\proce{\bar r_t^{\text{\tiny{D}}}}$ is the deterministic element of $\C\left([0,\infty),\R\right)$ 
defined for all $t$ by:
\begin{equation}
\label{eq:defovedrtD}
\bar r_t^{\text{\tiny{D}}}=\left\{d-t\car_{\left\{t \le \mu^{-1}\right\}}-\left(\frac{\rho-1}{\rho}t+\frac{1}{\lambda}\right)\car_{\left\{t > \mu^{-1}\right\}}\right\}^+.
\end{equation}
\end{proposition}
\begin{proof}
Fix $\xi>0$, assuming without loss of generality that: 
\begin{equation}
\label{eq:hypoxi0}
\xi < \frac{d-\mu^{-1}}{\rho}\wedge\frac{1}{\mu},
\end{equation} 
which implies that:
\begin{equation}
\label{eq:hypoxi1}
\mu^{-1}+\rho\xi<d<\bar \omega_0^{*,\text{\tiny{D}}}-\xi.
\end{equation}
If $\rho<3$, assume in addition to (\ref{eq:hypoxi0}) that:
\begin{equation}
\label{eq:hypoxi2}
\xi < \frac{1}{(4-\rho)\mu}.
\end{equation}
Let us first focus on
the interval of time $[0,\bar \omega_0^{*,\text{\tiny{D}}}-\xi].$ We have:
\begin{multline}
\label{eq:bach4}
\pr{\underset{0 \le t \le \bar \omega_0^{*,\text{\tiny{D}}}-\xi}{\sup} \left|t_1(\bar \nu_t^{n,\text{\tiny{D}}})-\bar r_t^{\text{\tiny{D}}}\right| >\xi}\\
 \shoveleft{\le \pr{\underset{0 \le t \le
       \left\{(\bar \omega_0^{*,\text{\tiny{D}}}-\xi)\wedge\bar \omega_0^{n,\text{\tiny{D}}}\wedge\bar \tau_0^{n,\text{\tiny{D}}}\right\}}{\sup}
     \left|t_1(\nu_{nt}^{\text{\tiny{D}}})-n\bar r_t^{\text{\tiny{D}}}\right| >n\xi}}\\
+\pr{\bar \omega_0^{*,\text{\tiny{D}}}-\bar \omega_0^{n,\text{\tiny{D}}} \ge \xi}+\pr{\tau_0^{n,\text{\tiny{D}}} \le
  n\bar \omega_0^{*,\text{\tiny{D}}}-n\xi}.
\end{multline}
Let us denote for all $y \in \R$ and $t \ge 0$:
\begin{equation*}
\tilde{\mathcal{ N}}^n_{t,y}:=\mbox{number of customers
    arrived up to time $t$, having deadline
    before $t+y$}.
\end{equation*} 
On the one hand, on the event 
\begin{equation*}
\left\{\underset{0 \le t \le
       \left\{(\bar \omega_0^{*,\text{\tiny{D}}}-\xi)\wedge\bar \omega_0^{n,\text{\tiny{D}}}\wedge\bar \tau_0^{n,\text{\tiny{D}}}\right\}}{\sup}
     \left(t_1(\nu_{nt}^{\text{\tiny{D}}})-n\bar r_t^{\text{\tiny{D}}}\right) >n\xi\right\},
\end{equation*}
there exists $t \in \left[0,\left\{(\bar \omega_0^{*,\text{\tiny{D}}}-\xi)\wedge\bar \omega_0^{n,\text{\tiny{D}}}\wedge\bar \tau_0^{n,\text{\tiny{D}}}\right\}\right]$ such that $t_1(\nu_{nt}^{\text{\tiny{D}}})-n\bar r_t^{\text{\tiny{D}}}>n\xi.$ 
Hence, since there has been no loss 
until $nt$, the number of services initiated up to $nt$
({i.e.}, $\mathcal S^n_{nt}+1$, which equals $S^n_{nt}+1$ in 
distribution) is larger than
the number of customers initially present, or arrived up to time $nt$, having deadline before
$nt+n\bar r_{t}^{\text{\tiny{D}}}+n\xi$ ({i.e.}, $\tilde{\mathcal{N}}^n_{nt,n\left(\bar r_t^{\text{\tiny{D}}}+\xi\right)}+n+1$).\\
On the other hand, on the event
\begin{equation}
\label{eq:event}
\left\{\underset{0 \le t \le
       \left\{(\bar \omega_0^{*,\text{\tiny{D}}}-\xi)\wedge\bar \omega_0^{n,\text{\tiny{D}}}\wedge\bar \tau_0^{n,\text{\tiny{D}}}\right\}}{\sup}
     \left(n\bar r_t^{\text{\tiny{D}}}-t_1(\nu_{nt}^{\text{\tiny{D}}})\right) >n\xi\right\},
\end{equation}
there exists $t \in \left[0,(\bar \omega_0^{*,\text{\tiny{D}}}-\xi)\wedge\bar \omega_0^{n,\text{\tiny{D}}}\wedge\bar \tau_0^{n,\text{\tiny{D}}}\right]$ such that $t_1(\nu_{nt}^{\text{\tiny{D}}})<n\bar r_t^{\text{\tiny{D}}}-n\xi.$ But for all $t \le (\mu^{-1})+\rho\xi$, every customer initially in the 
system, or arrived before $nt$ has a deadline at, or posterior to, $nd$, and 
hence a 
residual time credit at $nt$ larger or equal to $nd-nt \ge n\bar r_t^{\text{\tiny{D}}}-n\xi$. 
 Therefore, for the event (\ref{eq:event}) to 
occur, there must exist an instant 
$t \in \left[\mu^{-1}+\rho\xi,(\bar \omega_0^{*,\text{\tiny{D}}}-\xi)\wedge\bar \omega_0^{n,\text{\tiny{D}}}\wedge\bar \tau_0^{n,\text{\tiny{D}}}\right]$ such that $t_1(\nu_{nt}^{\text{\tiny{D}}})<n\bar r_t^{\text{\tiny{D}}}-n\xi.$ 
Since there has been no idle time on the interval 
$[0,\omega_0^{n,\text{\tiny{D}}}\wedge\tau_0^{n,\text{\tiny{D}}}]$, and since the discipline amounts to FIFO, this 
implies that $\mathcal S^n_{nt}+1$ (equal in distribution to $S^n_{nt}+1$) is 
less than the number of customers arrived up to time $nt$, having
deadline before $nt+n\bar r_{t}^{\text{\tiny{D}}}-n\xi$ (that is, $\tilde{\mathcal{N}}^n_{nt,n\left(\bar r_t^{\text{\tiny{D}}}-\xi\right)}+n+1$). 
Consequently:
\begin{multline}
\label{eq:bach5}
\pr{\underset{0 \le t \le
       \left\{(\bar \omega_0^{*,\text{\tiny{D}}}-\xi)\wedge\bar \omega_0^{n,\text{\tiny{D}}}\wedge\bar \tau_0^{n,\text{\tiny{D}}}\right\}}{\sup}
     \left|t_1(\nu_{nt}^{\text{\tiny{D}}})-n\bar r_t^{\text{\tiny{D}}}\right| > n\xi}\\
\shoveleft{\le \pr{\underset{0 \le t \le
       \left\{(\bar \omega_0^{*,\text{\tiny{D}}}-\xi)\wedge\bar \omega_0^{n,\text{\tiny{D}}}\wedge\bar \tau_0^{n,\text{\tiny{D}}}\right\}}{\sup}
     \left(t_1(\nu_{nt}^{\text{\tiny{D}}})-n\bar r_t^{\text{\tiny{D}}}\right) >n\xi}}\\
\shoveright{+\pr{\underset{0 \le t \le
       \left\{(\bar \omega_0^{*,\text{\tiny{D}}}-\xi)\wedge\bar \omega_0^{n,\text{\tiny{D}}}\wedge\bar \tau_0^{n,\text{\tiny{D}}}\right\}}{\sup}
     \left(n\bar r_t^{\text{\tiny{D}}}-t_1(\nu_{nt}^{\text{\tiny{D}}})\right) >n\xi}}\\
\shoveleft{\le \pr{\underset{0 \le t \le
       \bar \omega_0^{*,\text{\tiny{D}}}-\xi}{\sup}
     \left(S^n_{nt}-\tilde{\mathcal{N}}^n_{nt,n(\bar r_t^{\text{\tiny{D}}}+\xi)}\right)\ge n
   }}\\+\pr{\underset{\mu^{-1}+\rho\xi \le t \le
       \bar \omega_0^{*,\text{\tiny{D}}}-\xi}{\sup}
     \left(\tilde{\mathcal{N}}^n_{nt,n(\bar r_t^{\text{\tiny{D}}}-\xi)}-S^n_{nt}\right)\ge -n}.
\end{multline}
It is easily checked, that for all $n$ and
all $z \in \R$,   
$$\proce{\tilde {\mathcal {N}}^n_{nt,n(\bar r_t^{\text{\tiny{D}}}+z)}-\tilde \lambda^n_{t,z}}$$ is a $\mathcal
G^n_t$-martingale, where for all $t\ge 0$:
\begin{equation*}
\tilde\lambda^n_{t,z}=n\lan\left\{t-\left(d-(\bar r_t^{\text{\tiny{D}}}+z)\right)^+\right\}^+.
\end{equation*}
Thus, the processes  
$$\tilde{\mathcal
  {M}}^n_{\xi}(t):=S^n_{nt}-\tilde{\mathcal{N}}^n_{nt,n\left(\bar r_t^{\text{\tiny{D}}}+\xi\right)}-\left(\mu^nnt-\tilde
\lambda^n_{t,\xi}\right)$$
and 
\begin{equation*}
\widehat{\mathcal
  {M}}^n_{\xi}(t):=\tilde{\mathcal{N}}^n_{nt,n\left(\bar r_t^{\text{\tiny{D}}}-\xi\right)}-S^n_{nt}-\left(\tilde
\lambda^n_{t,-\xi}-\mu^nnt\right)
\end{equation*}
are $\mathcal G^n_t-$
martingales. 
Hence, (\ref{eq:bach5}) becomes:
\begin{multline}
\label{eq:bach5bis}
\pr{\underset{0 \le t \le
       \left\{(\bar \omega_0^{*,\text{\tiny{D}}}-\xi)\wedge\bar \omega_0^{n,\text{\tiny{D}}}\wedge\bar \tau_0^{n,\text{\tiny{D}}}\right\}}{\sup}
     \left|t_1(\nu_{nt}^{\text{\tiny{D}}})-n\bar r_t^{\text{\tiny{D}}}\right| >n\xi}\\ 
\shoveleft{\le \pr{\underset{0 \le t \le
       \bar \omega_0^{*,\text{\tiny{D}}}-\xi}{\sup}
     \tilde{\mathcal {M}}^n_{\xi}(t) \ge n+\underset{0 \le t \le
       \bar \omega_0^{*,\text{\tiny{D}}}-\xi}{\inf}\left(\tilde
     \lambda^n_{t,\xi}-\mu^nnt\right) }}\\+\pr{\underset{0 \le t \le
       \bar \omega_0^{*,\text{\tiny{D}}}-\xi}{\sup}
     \widehat{\mathcal {M}}^n_{\xi}(t)\ge -n+\underset{\mu^{-1}+\rho\xi \le t \le
       \bar \omega_0^{*,\text{\tiny{D}}}-\xi}{\inf}\left(\mu^nnt-\tilde
     \lambda^n_{t,-\xi}\right)}.
\end{multline}
On the one hand, since:
\[\tilde
     \lambda^n_{t,\xi}-\mu^nnt=\left\{\begin{array}{ll}
 n\left\{\lambda^n(t\wedge \xi)-\mu^nt\right\} &\mbox{ if }t\in \left[0,\mu^{-1}\right]\\
 n\left\{\frac{\lambda^n\mu}{\lambda}t-\frac{\lambda^n}{\lambda}+\lambda\n\xi-\mu^n t \right\}
&\mbox{ if }t \in \left[\mu^{-1},\bar \omega_0^{*,\text{\tiny{D}}}-\xi\right], 
\end{array}
\right.\] 
it follows from Hypothesis \ref{hypo:HDET} that for a sufficiently large $n$, 
 $$n+\underset{0 \le t \le
       \bar \omega_0^{*,\text{\tiny{D}}}-\xi}{\inf}\left(\tilde
     \lambda^n_{t,\xi}-\mu^nnt\right) \ge n+n\left\{\frac{\lambda\xi}{2}-1\right\}>0.$$
On the other hand, since for all $t \in \left[\mu^{-1}+\rho\xi,\bar \omega_0^{*,\text{\tiny{D}}}-\xi\right]$, 
$$\mu^nnt-\tilde
     \lambda^n_{t,-\xi}=n\left\{\mu^nt-\frac{\lambda^n\mu}{\lambda}t+\frac{\lambda^n}{\lambda}+\lambda^n\xi\right\},$$
for a sufficiently large $n$,  
$$-n+\underset{\mu^{-1}+\rho\xi \le t \le
       \bar \omega_0^{*,\text{\tiny{D}}}-\xi}{\inf}\left(\mu^nnt-\tilde
     \lambda^n_{t,-\xi}\right)\ge -n +n\left\{\frac{\lambda\xi}{2}+1\right\}>0.$$ 
Thus from (\ref{eq:bach5bis}), for a sufficiently large $n$, 
\begin{multline}
\label{eq:bach6}
\pr{\underset{0 \le t \le
       \left\{(\bar \omega_0^{*,\text{\tiny{D}}}-\xi)\wedge\bar \omega_0^{n,\text{\tiny{D}}}\wedge\bar \tau_0^{n,\text{\tiny{D}}}\right\}}{\sup}
     \left|t_1(\nu_{nt}^{\text{\tiny{D}}})-n\bar r_t^{\text{\tiny{D}}}\right| > n\xi}\\
\shoveleft{\le \frac{16}{n\lambda\xi^2}\left\{\esp{\left\{\tilde{\mathcal {M}}^n_{\xi}\left(\bar \omega_0^{*,\text{\tiny{D}}}-\xi\right)\right\}^2}+\esp{\left\{\widehat{\mathcal {M}}^n_{\xi}\left(\bar \omega_0^{*,\text{\tiny{D}}}-\xi\right)\right\}^2}\right\}}\\
= \frac{16}{\left(n\lambda\xi\right)^2}\left\{\tilde \lambda^n_{(\bar \omega_0^{*,\text{\tiny{D}}}-\xi),\xi}+\tilde
  \lambda^n_{(\bar \omega_0^{*,\text{\tiny{D}}}-\xi),-\xi}+2\mu^n
  n(\bar \omega_0^{*,\text{\tiny{D}}}-\xi)\right\}\tend 0,
\end{multline}
using successively Tchebitchef and Doob's inequalities. \\\\
Consider now the term: 
\begin{multline}
\label{eq:bach7}
\pr{\tau_0^{n,\text{\tiny{D}}} \le
  n\bar \omega_0^{*,\text{\tiny{D}}}-n\xi}
\le \pr{\omega_0^{n,\text{\tiny{D}}} \le n\bar \omega_0^{*,\text{\tiny{D}}}-n\xi}\\+\pr{\left\{\tau_0^{n,\text{\tiny{D}}}\le n\bar \omega_0^{*,\text{\tiny{D}}}-n\xi\right\}\cap\left\{\omega_0^{n,\text{\tiny{D}}}>n\bar \omega_0^{*,\text{\tiny{D}}}-n\xi\right\}}.
\end{multline}
On the event $$\left\{\tau_0^{n,\text{\tiny{D}}}\le
  n\bar \omega_0^{*,\text{\tiny{D}}}-n\xi\right\}\cap\left\{\omega_0^{n,\text{\tiny{D}}}>n\bar \omega_0^{*,\text{\tiny{D}}}-n\xi\right\},$$
there exists an instant, say $t \in [0,\bar \omega_0^{*,\text{\tiny{D}}}-\xi]$ such that
the buffer is empty at $nt$, but there has been no loss 
before $nt$. Applying the same arguments that led to (\ref{eq:bach3})
yields: 
\begin{equation*}
\pr{\left\{\tau_0^{n,\text{\tiny{D}}}\le
    n\bar \omega_0^{*,\text{\tiny{D}}}-n\xi\right\}\cap\left\{\omega_0^{n,\text{\tiny{D}}}>n\bar \omega_0^{*,\text{\tiny{D}}}-n\xi\right\}}\tend 0,
\end{equation*}
which together with \lemref{lemma:convomega} and (\ref{eq:bach7}) yields:
\begin{equation}
\label{eq:bach8}
\pr{\tau_0^{n,\text{\tiny{D}}} \le
  n\bar \omega_0^{*,\text{\tiny{D}}}-n\xi}\tend 0.
\end{equation}
Finally, using (\ref{eq:bach6}) together with (\ref{eq:bach8}) and
\lemref{lemma:convomega} in (\ref{eq:bach4}):
\begin{equation}
\label{eq:bach9}
\pr{\underset{0 \le t \le \bar \omega_0^{*,\text{\tiny{D}}}-\xi}{\sup} \left|t_1(\bar \nu_t^{n,\text{\tiny{D}}})-\bar r_t^{\text{\tiny{D}}}\right| >
  \xi}\tend 0.
\end{equation}
\textit{ }\\\\\\
Let us now consider the interval of time
$[\bar \omega_0^{*,\text{\tiny{D}}}-\xi,\infty).$ 
We have:
\begin{multline}
\label{eq:ange0}
\pr{\underset{t \ge \bar \omega_0^{*,\text{\tiny{D}}}-\xi}{\sup} \left|t_1(\bar \nu_t^{n,\text{\tiny{D}}})-\bar r_t^{\text{\tiny{D}}}\right| > 2\xi}\\
\shoveleft{\le \pr{\underset{t \ge \bar \omega_0^{*,\text{\tiny{D}}}-\xi }{\sup} \left(t_1(\nu_{nt}^{\text{\tiny{D}}})-n\bar r_t^{\text{\tiny{D}}}\right) > 2n\xi}}\\
+\pr{\underset{t \ge \bar \omega_0^{*,\text{\tiny{D}}}-\xi 
      }{\inf} t_1(\nu_{nt}^{\text{\tiny{D}}}) < n\left(\underset{t \ge \bar \omega_0^{*,\text{\tiny{D}}}-\xi }{\sup}\bar r_t^{\text{\tiny{D}}}-2\xi\right)}.
\end{multline}
First, it is easily seen that $$\underset{t \ge \bar \omega_0^{*,\text{\tiny{D}}}-\xi }{\sup}\bar r_t^{\text{\tiny{D}}}=\bar r^{\text{\tiny{D}}}_{\bar \omega_0^{*,\text{\tiny{D}}}-\xi}=\frac{\rho-1}{\rho}\xi<\xi.$$ 
Therefore:
\begin{equation}
\label{eq:ange1}
\pr{\underset{t \ge \bar \omega_0^{*,\text{\tiny{D}}}-\xi }{\inf} t_1(\nu_{nt}^{\text{\tiny{D}}}) < n\left(\underset{t \ge \bar \omega_0^{*,\text{\tiny{D}}}-\xi }{\sup}\bar r_t^{\text{\tiny{D}}}-2\xi\right)}=0.
\end{equation}
\textit{ }\\\\
On another hand, 
define the following event: 
\begin{multline*}
\mathcal E^n_{\xi}:=\Biggl\{\text{for all $t\ge \bar \omega_0^{*,\text{\tiny{D}}}-\xi$, some customers arrive before
      $n(t-\xi)$,}\\\shoveright{\text{ with deadline in $\left[n(t+\bar r_t^{\text{\tiny{D}}}+\xi), n(t+\bar r_t^{\text{\tiny{D}}}+2\xi)\right]$}\Biggl\}}\\
=\left\{\underset{t\ge \bar \omega_0^{*,\text{\tiny{D}}}-\xi}{\inf}\left\{\tilde{\mathcal{N}}^n_{n(t-\xi),n\left(\bar r_{t-\xi}^{\text{\tiny{D}}}+\frac{2\rho+1}{\rho}\xi\right)}-\tilde{\mathcal{N}}^n_{n(t-\xi),n\left(\bar r_{t-\xi}^{\text{\tiny{D}}}+\frac{\rho+1}{\rho}\xi\right)}\right\}>0\right\}.
\end{multline*}
We have: 
\begin{multline}
\label{eq:ange31}
\pr{\underset{t \ge \bar \omega_0^{*,\text{\tiny{D}}}-\xi }{\sup}
  \left(t_1(\nu_{nt}^{\text{\tiny{D}}})-\bar r_t^{\text{\tiny{D}}}\right) > 2n\xi}\\
\shoveleft{\le\mathbf P\Biggl[\left\{\underset{t \ge \bar \omega_0^{*,\text{\tiny{D}}}-\xi }{\sup}
  \left(t_1(\nu_{nt}^{\text{\tiny{D}}}) -n\bar r_t^{\text{\tiny{D}}}\right)> 2n\xi\right\}}\\\shoveright{\cap\left\{\underset{0\le
    t \le \bar \omega_0^{*,\text{\tiny{D}}}-\xi}{\sup}\left|t_1(\nu_{nt}^{\text{\tiny{D}}})-n\bar r_t^{\text{\tiny{D}}}\right| \le 
  n\xi\right\}\cap \mathcal E^n_{\xi}\Biggl]}\\
+\pr{\underset{0\le
    t \le \bar \omega_0^{*,\text{\tiny{D}}}-\xi}{\sup}\left|t_1(\bar \nu_t^{n,\text{\tiny{D}}})-\bar r_t^{\text{\tiny{D}}}\right|>
  \xi}+\pr{\left(\mathcal E^n_{\xi}\right)^c}.
\end{multline}
On the event \begin{multline*}\left\{\underset{t \ge \bar \omega_0^{*,\text{\tiny{D}}}-\xi }{\sup}
  \left(t_1(\nu_{nt}^{\text{\tiny{D}}}) -n\bar r_t^{\text{\tiny{D}}}\right)> 2n\xi\right\}\\\cap\left\{\underset{0\le
    t \le \bar \omega_0^{*,\text{\tiny{D}}}-\xi}{\sup}\left|t_1(\nu_{nt}^{\text{\tiny{D}}})-n\bar r_t^{\text{\tiny{D}}}\right|\le
  n\xi\right\}\cap \mathcal E^n_{\xi} ,
\end{multline*} 
there exists an instant $t \ge \bar \omega_0^{*,\text{\tiny{D}}}-\xi$ (the
first one), such
that there is no customer in the system at $nt$ having deadline
between $n(t+\bar r_t^{\text{\tiny{D}}}+\xi)$ and $n(t+\bar r_t^{\text{\tiny{D}}}+2\xi).$ Consider the customer who 
would have had the smallest time credit at $nt$ if we was still present in 
the system at this instant, among those arrived 
before $n(t-\xi)$, with deadline in
$\left[n(t+\bar r_t^{\text{\tiny{D}}}+\xi),n(t+\bar r_t^{\text{\tiny{D}}}+2\xi)\right]$ 
. For all
$s\le t-\xi$ such that this customer has already entered the system at $ns$, denote by $\tilde R_s$, the
time credit of this customer at $ns$. We have: 
$$\tilde R_s \in \biggl[n(t-s)+n\bar r_t^{\text{\tiny{D}}}+n\xi,n(t-s)+n\bar r_t^{\text{\tiny{D}}}+2n\xi\biggl].$$
In particular, it is easily seen with the form of $\bar r_t^{\text{\tiny{D}}}$ that:
\begin{equation*}
\tilde R_s \ge n(t-s)+n\bar r_t^{\text{\tiny{D}}}+n\xi
> n\bar r_s^{\text{\tiny{D}}}+n\xi \ge t_1(\nu_{ns}^{n,\text{\tiny{D}}})
\end{equation*} 
(the last inequality is true since $s \le \bar \omega_0^{*,\text{\tiny{D}}}-\xi$ and $$\underset{0\le
    u \le \bar \omega_0^{*,\text{\tiny{D}}}-\xi}{\sup}\left|t_1(\nu_{nu}^{n,\text{\tiny{D}}})-n\bar r_u^{\text{\tiny{D}}}\right|\le
  n\xi\,\,).$$
Thus, none of the customers arrived before $n(t-\xi)$, with
deadline  in\\ 
$\left[n(t+\bar r_t^{\text{\tiny{D}}}+\xi),n(t+\bar r_t^{\text{\tiny{D}}}+2\xi)\right]$
have been served before $n(t-\xi)$, since none of them has ever
been prioritary. Hence:
\begin{multline}
\label{eq:ange33quart}
\mathbf P\Biggl[\left\{\underset{t \ge \bar \omega_0^{*,\text{\tiny{D}}}-\xi }{\sup}
  \left(t_1(\nu_{nt}^{\text{\tiny{D}}})-n\bar r_t^{\text{\tiny{D}}}\right) >2n\xi\right\}\\\shoveright{\cap\left\{\underset{\varepsilon\le
    t \le \bar \omega_0^{*,\text{\tiny{D}}}-\xi}{\sup}\left|t_1(\nu_{nt}^{\text{\tiny{D}}})-n\bar r_t^{\text{\tiny{D}}}\right|\le n\xi\right\}\cap \mathcal E^n_{\xi}\Biggl]}\\
\shoveleft{\le \mathbf P\Biggl[\biggl\{\mbox{All of the customers arrived before $n(t-\xi)$,}}\\\mbox{ with deadline in
$\left[n(t+\bar r_t^{\text{\tiny{D}}}+\xi),n(t+\bar r_t^{\text{\tiny{D}}}+2\xi)\right]$}\\\shoveright{\mbox{
  have entered service between $n(t-\xi)$ and $nt$}\biggl\}\cap \mathcal E^n_{\xi}\Biggl]}\\
\shoveleft{\le \mathbf P\Biggl[\underset{t \ge \bar \omega_0^{*,\text{\tiny{D}}}-\xi}{\sup}
  \Biggl\{S^n_{n\xi}-\biggl(\tilde{\mathcal{N}}^n_{n(t-\xi),n\left(\bar r_{t-\xi}^{\text{\tiny{D}}}+\frac{2\rho+1}{\rho}\xi\right)}}\\-\tilde{\mathcal{N}}^n_{n(t-\xi),n\left(\bar r_{t-\xi}^{\text{\tiny{D}}}+\frac{\rho+1}{\rho}\xi\right)}\biggl)\Biggl\}\ge 0\Biggl],
\end{multline} 
where $S^n_{n\xi}$ denotes the value at $n\xi$ of a Poisson process of
intensity $\mu^n$. The following is a $\mathcal G^n_t-$ martingale:
\begin{multline*}
\check{\mathcal{M}}^n_{t,\xi}:=
S^n_{n\xi}-\left\{\tilde{\mathcal{N}}^n_{n(t-\xi),n\left(\bar r_{t-\xi}^{\text{\tiny{D}}}+\frac{2\rho+1}{\rho}\xi\right)}-\tilde{\mathcal{N}}^n_{n(t-\xi),n\left(\bar r_{t-\xi}^{\text{\tiny{D}}}+\frac{\rho+1}{\rho}\xi\right)}\right\}\\-\left\{\mu^nn\xi-\left(\tilde
    \lambda^n_{(t-\xi),\frac{2\rho+1}{\rho}\xi}-\tilde
    \lambda^n_{(t-\xi),\frac{\rho+1}{\rho}\xi}\right)\right\},
\end{multline*} 
and for all $\xi$, $\check{\mathcal{M}}^n_{t,\xi} \in L^2$. 
Then, it is easily seen that  the last term of (\ref{eq:ange33quart}) is equal to:
\begin{equation*}
\pr{\underset{t \ge \bar \omega_0^{*,\text{\tiny{D}}}-\xi}{\sup}\check{\mathcal{M}}^n_{t,\xi}
  \ge \underset{t \ge \bar \omega_0^{*,\text{\tiny{D}}}-\xi}{\inf}\left\{\left(\tilde
    \lambda^n_{(t-\xi),\frac{2\rho+1}{\rho}\xi}-\tilde
    \lambda^n_{(t-\xi),\frac{\rho+1}{\rho}\xi}\right)-\mu^nn\xi\right\}}.
\end{equation*}
But for all $t\ge\bar \omega_0^{*,\text{\tiny{D}}}-\xi$:
\begin{multline*}
\tilde
    \lambda^n_{(t-\xi),\frac{2\rho+1}{\rho}\xi}-\tilde
    \lambda^n_{(t-\xi),\frac{\rho+1}{\rho}\xi}\\
\shoveleft{=\lambda^n\left\{t-\xi-\left(d-\bar r_{t-\xi}^{\text{\tiny{D}}}-\frac{2\rho+1}{\rho}\xi\right)^+\right\}^+}\\-\lambda^n\left\{t-\xi-\left(d-\bar r_{t-\xi}^{\text{\tiny{D}}}-\frac{\rho+1}{\rho}\xi\right)^+\right\}^+
=\lambda^n \xi,
\end{multline*}
using (\ref{eq:hypoxi1}) and (\ref{eq:hypoxi2}).  
Consequently, for a sufficiently large $n$:
\begin{multline}
\label{eq:ange33}
\mathbf P\Biggl[\left\{\underset{t \ge \bar \omega_0^{*,\text{\tiny{D}}}-\xi }{\sup}
  \left(t_1(\nu_{nt}^{\text{\tiny{D}}})-n\bar r_t^{\text{\tiny{D}}}\right) >2n\xi\right\}\\\shoveright{\cap\left\{\underset{\varepsilon\le
    t \le \bar \omega_0^{*,\text{\tiny{D}}}-\xi}{\sup}\left|t_1(\nu_{nt}^{\text{\tiny{D}}})-n\bar r_t^{\text{\tiny{D}}}\right|\le n\xi\right\}\cap \mathcal E^n_{\xi}\Biggl]}\\
\shoveleft{\le \pr{\underset{t \ge \bar \omega_0^{*,\text{\tiny{D}}}-\xi}{\sup}\check{\mathcal{M}}^n_{t,\xi}
  \ge n\left(\lambda^n-\mu^n\right)\xi}}\\
\le\frac{4}{\left(n\left(\lambda^n-\mu^n\right)\xi\right)^2}\left\{\lambda^nn\xi+\mu^nn\xi\right\}\tend 0,
\end{multline}
using again Doob inequality. 
Remark, that we also proved that:
\begin{multline}
\label{eq:ange34}
\pr{\left(\mathcal E^n_{\xi}\right)^c}\\=
\pr{\underset{t \ge \bar \omega_0^{*,\text{\tiny{D}}}-\xi}{\inf}
 \left(\tilde{\mathcal{N}}^n_{n(t-\xi),n\left(\bar r_{t-\xi}^{\text{\tiny{D}}}+\frac{2\rho+1}{\rho}\xi\right)}-\tilde{\mathcal{N}}^n_{n(t-\xi),n\left(\bar r_{t-\xi}^{\text{\tiny{D}}}+\frac{\rho+1}{\rho}\xi\right)}\right\}\le 0}\\
\shoveleft{\le \mathbf P\Biggl[\underset{t \ge \bar \omega_0^{*,\text{\tiny{D}}}-\xi}{\sup}
  \Biggl\{S^n_{n\xi}-\biggl(\tilde{\mathcal{N}}^n_{n(t-\xi),n\left(\bar r_{t-\xi}^{\text{\tiny{D}}}+\frac{2\rho+1}{\rho}\xi\right)}}\\\shoveright{-\tilde{\mathcal{N}}^n_{n(t-\xi),n\left(\bar r_{t-\xi}^{\text{\tiny{D}}}+\frac{\rho+1}{\rho}\xi\right)}\biggl)\Biggl\}\ge 0\Biggl]}\\ \tend 0.
\end{multline} 
Therefore, using (\ref{eq:ange33}), (\ref{eq:ange34}) and (\ref{eq:bach9}) in
(\ref{eq:ange31}) yields to:
\begin{equation*}
\pr{\underset{t \ge \bar \omega_0^{*,\text{\tiny{D}}}-\xi }{\sup}
  \left(t_1(\nu_{nt}^{\text{\tiny{D}}})-\bar r_t^{\text{\tiny{D}}}\right) > 2n\xi}
 \tend 0.
\end{equation*}
Together with (\ref{eq:ange1}) in (\ref{eq:ange0}), this entails:
\begin{equation*}
\pr{\underset{t \ge \bar \omega_0^{*,\text{\tiny{D}}}-\xi}{\sup} \left|t_1(\bar \nu_t^{n,\text{\tiny{D}}})-\bar r_t^{\text{\tiny{D}}}\right| > 2\xi}\tend 0.
\end{equation*}
This implies in turns, together with (\ref{eq:bach9}), that:
\begin{multline*}
\pr{\underset{t \ge
    0}{\sup}\left|t_1(\bar \nu_t^{n,\text{\tiny{D}}})-\bar r_t^{\text{\tiny{D}}}\right|\ge 3\xi}
\le \pr{\underset{0 \le t \le \bar \omega_0^{*,\text{\tiny{D}}}-\xi}{\sup} \left|t_1(\bar \nu_t^{n,\text{\tiny{D}}})-\bar r_t^{\text{\tiny{D}}}\right| > \xi}\\\shoveright{+\pr{\underset{t \ge \bar \omega_0^{*,\text{\tiny{D}}}-\xi}{\sup} \left|t_1(\bar \nu_t^{n,\text{\tiny{D}}})-\bar r_t^{\text{\tiny{D}}}\right| > 2\xi}}\\ \tend 0.
\end{multline*}
We conclude using ~\cite{Rob}.
\end{proof}
We can therefore conclude with the next result, which yields then
convergence in distribution of a normalized sequence of profile
processes of a M/M/1/1+D-EDF queue to an explicit fluid limit:

\begin{theorem}[Fluid limit of the M/M/1/1+D-EDF queue]
\label{thm:ouiDET}
{
$$\suite{\proce{\bar \nu_t^{n,\text{\tiny{D}}}}} \Longrightarrow \proce{\bar \nu_t^{*,\text{\tiny{D}}}}\mbox{ , in $\D\left([0,\infty),\M_f^+\right)$},$$
where for all $t \ge 0$ and all $f \in \D_b$: 
\begin{equation*}
\cro{\bar \nu^{*,\text{\tiny{D}}}_t,f}=f(d-t)-\mu\int_0^t f\left(\bar
  r^{\text{\tiny{D}}}_{s}+s-t\right)\,ds
+\lambda \int_0^t f\left(d-s\right)\,ds,
\end{equation*}
$\proce{\bar r_t^{\text{\tiny{D}}}}$ being defined by (\ref{eq:defovedrtD}).}
\end{theorem}
\begin{proof}
Let us verify the assumptions of Theorem \ref{thm:oui} for any given $T \ge 0$. First, it is straightforward, that Hypothesis \ref{hypo:H1} are satisfied in view of Hypothesis \ref{hypo:HDET}, for $\bar D=d$, a.s., and for 
$\xi$ the Dirac measure at $d$. 

Then, remark, that in (\ref{eq:ange33}) we proved as a matter of fact
that for all $\xi>0$ satisfying (\ref{eq:hypoxi0}) and (\ref{eq:hypoxi2}), 
\begin{multline*}
\mathbf P\Biggl[\biggl\{\mbox{There exists $t \ge
      \bar \omega_0^{*,\text{\tiny{D}}}-\xi$, such that all of the customers arrived 
    }\\\mbox{before $n(t-\xi)$, with deadline in
$\left[n(t+\bar r_t^{\text{\tiny{D}}}+\xi),n(t+\bar r_t^{\text{\tiny{D}}}+2\xi)\right]$}\\\shoveright{\mbox{
  have been served before $nt$}\biggl\}\cap\mathcal
E^n_{\xi}\Biggl]}\\\tend 0.
\end{multline*} Denoting this previous event $\mathcal A^n_{\xi}$, 
it follows that:
\begin{equation*}
\pr{\bar \tau_0^{n,\text{\tiny{D}}} \le T}\le \pr{\bar \tau_0^{n,\text{\tiny{D}}} \le \bar \omega_0^{*,\text{\tiny{D}}}-\xi}+\pr{\mathcal A^n_{\xi}}+\pr{\left(\mathcal
E^n_{\xi}\right)^c}\tend 0,
\end{equation*}
using (\ref{eq:ange34}) and (\ref{eq:bach8}) as well: (\ref{eq:TpsVidage}) is satisfied. 
Finally, (\ref{eq:ConvPre}) is verified in view of Proposition \ref{pro:ConvPre} for $\proce{\bar r_t^{\text{\tiny{D}}}}$ defined in (\ref{eq:defovedrtD}). 
We therefore can apply Theorem \ref{thm:oui} for all $T\ge 0$, which completes the proof.
\end{proof}
We can in particular approximate the queue length and loss processes, by applying the profile process of the queue to simple rcll functions. First, 
in view of (\ref{eq:normQ}), the normalized queue length process 
$\proc{\bar Q^n}$ can 
be asymptotically approximated by the process defined for all $t \ge 0$ by:
\begin{multline*}
\cro{\bar \nu_t^{*,\text{\tiny{D}}},\car_{\Rp}}=-\mu\int_0^t \car_{\Rp}\left(\bar
  r^{\text{\tiny{D}}}_{s}+s-t\right)\,ds
+\lambda \int_0^t \car_{\Rp}\left(d-s\right)\,ds\\
=\left\{1+(\lambda-\mu)t\right\}\car_{\left\{t \le \frac{\rho d-\mu^{-1}}{\rho -
      1}\right\}}+\lambda d \car_{\left\{t \ge \frac{\rho d-\mu^{-1}}{\rho -
      1}\right\}}.
\end{multline*}
Similarly, in view of (\ref{eq:normP}) the normalized loss process 
$\proc{\bar P^n}$ can be approximated 
 by the process defined for all $t \ge 0$ by:
\begin{multline*}
\cro{\bar \nu_t^{*,\text{\tiny{D}}},\car_{\R-}}=-\mu\int_0^t \car_{\R-}\left(\bar
  r^{\text{\tiny{D}}}_{s}+s-t\right)\,ds
+\lambda \int_0^t \car_{\R-}\left(d-s\right)\,ds\\
=\left(1+\lambda(t-d)-\mu t\right) \car_{\left\{t \ge \frac{\rho
      d-\mu^{-1}}{\rho -
      1}\right\}}.
\end{multline*}

\subsection*{{M/GI/$\infty$ system}}
Remark, that we can also apply \thmref{thm:oui} to obtain the fluid limit of a pure delay M/GI/$\infty$ system. 
For all $n \in \N^*$, consider a {pure delay} (PD) M/GI/$\infty$ system: customers arrive according to 
a Poisson process of intensity $\lambda^n$, requesting service durations,  
i.i.d. of the distribution of $\alpha^n$,  
to an infinite reservoir of servers, assuming that:
$$\lambda^n \rightarrow \lambda,\,\esp{\frac{\alpha^n}{n}}\rightarrow \esp{\alpha},\,\frac{\alpha^n}{n} \overset{\mathcal D}{\rightarrow} \alpha.$$     
Each customer is hence immediately attended upon arrival, and remains in the system for the duration of his service. 
Such a system can be described by a profile process, keeping now track of the remaining processing 
times of the customers in service (positive atoms) and the elapsed times since departure of the already served customers 
(negative atoms):  
this is the profile process $\proce{\nu_t^{n,\text{\tiny{PD}}}}$ associated to 
$\left(\nu_0^{n,\text{\tiny{PD}}},\lambda^n,0,\alpha^n\right)$, where $\nu_0^{n,\text{\tiny{PD}}}$ is the 
profile of the service durations of the customers initially in the system 
(which we assume to be an $n$-sample of the distribution of $\alpha^n$). The service durations replace the time credits 
of the queue with impatient customers, and consequently the analogous service rate becomes null in this case. 
We normalize this process as above, writing for all $t\ge 0$ and all Borel 
set $\B$:
$$\bar \nu_t^{n,\text{\tiny{PD}}}(\B)=\frac{\nu_{nt}^{n,\text{\tiny{PD}}}(n\B)}{n}.$$ 
The fact that $\mu^n$ is zero allows us to skip conditions (\ref{eq:TpsVidage}) and (\ref{eq:ConvPre}) in \thmref{thm:oui}, whose application becomes 
straightforward: 
$$\suite{\proce{\ovnunPD_t}}\Longrightarrow \proce{\ovnuPD_t}\mbox{ in }\D\left([0,\infty),\M_f^+\right),$$ 
 where $\proc{\ovnuPD}$ is the deterministic element of
 $\D_{\infty,\M_f^+}$ defined for all $t\ge 0$ and all $f \in \D_b$ by:
   \begin{equation*}
     \cro{\ovnuPD_t,f}=\esp{\tau_tf(\alpha)}+\lambda\int_0^t\esp{\tau_sf(\alpha)}\,ds.
   \end{equation*} 
Hence, as above, we can asymptotically approximate the normalized congestion 
process (number of customers in service), given for all $n \in \N^*$ by 
$\proce{\cro{\ovnunPD_t,\car_{\Rp}}}$ by the process defined for all 
$t\ge 0$ by:
\begin{equation*}
\cro{\ovnuPD_t,\car_{\Rp}}=\pr{\alpha > t}
+\lambda\int_0^t\pr{\alpha > s}\,ds. 
\end{equation*} 
Remark, that in the special case where $\alpha^n$ is exponentially distributed 
of parameter $\mu^n \rightarrow \mu$, this process becomes for all $t\ge 0:$ 
$$e^{-\mu t}+\rho\left(1-e^{-\mu t}\right),$$
which is the fluid limit obtained by Borovkov in ~\cite{Bor67}.
We can also approximate the normalized workload 
process, given for all $n$ by $\proce{\cro{\ovnunPD_t,I\car_{\Rp}}}$ by the process given by:
\begin{equation*}
\cro{\ovnuPD_t,I\car_{\Rp}}=\esp{(\alpha -t)^+}+\lambda\int_0^t\esp{(\alpha -s)^+}\,ds. 
\end{equation*} 
Finally, the normalized process counting the already served customers can be approximated by the process defined for all $t\ge 0$ by:
$$\cro{\ovnuPD_t,\car_{\R-}}=\pr{\alpha \le t}
+\lambda\int_0^t\pr{\alpha \le s}\,ds.$$

\bibliographystyle{amsalpha}
\bibliography{bibliographie4}

\end{document}